# Learning Complexity Dimensions for a Continuous-Time Control System

Pirkko Kuusela, Daniel Ocone, and Eduardo D. Sontag
Rutgers, The State University of New Jersey, USA

October 24, 2018

## Abstract

This paper takes a computational learning theory approach to a problem of linear systems identification. It is assumed that input signals have only a finite number $k$ of frequency components, and systems to be identified have dimension no greater than $n$. The main result establishes that the sample complexity needed for identification scales polynomially with $n$ and logarithmically with $k$.

## 1 Introduction

The problem of systems identification may be seen as an instance of the general question of "learning" an unknown function. Thus, one may ask if techniques from Computational Learning Theory ("CLT" in what follows; see for instance [20] for a general introduction, as well as the papers in the Special Issue [4] of Systems and Control Letters) can be used to obtain insight into this central question in control theory.

Indeed, Ljung asked precisely this question in a paper presented at a Special Session in Learning Theory organized at the 1996 Conference on Decision and Control, see [14], and, independently, the papers [7, 8] had already provided results along these lines for discrete-time linear systems on finite-window data. See also [11], [9], and references there, as well as [12] for several results for nonlinear systems in discrete time.

For continuous-time systems, the situation is complicated by the fact that, even for finite-length inputs, learnability is impossible when formulated in the CLT framework; this was proved in [23], and, alternatively, can be seen by applying the discrete-time results (through sampling) from [7, 8].

Thus, in this paper, we suppose that all inputs to be used, in the learning as well as validation stages, belong to the linear span of a fixed number $k$ of sinusoidal basic functions. This band-limiting assumption allows us to obtain a precise result: the sample complexity needed for identification scales polynomially on an upper bound on the systems being identified, and logarithmically with $k$. This provides a tight analogy to the discrete results previously obtained, in which $k$ appeared as the length of the discrete-time window employed.

The reader is referred to [17] for results that apply to a class of nonlinear continuous-time systems, but which is formulated in terms of learning derivatives evaluated at a particular instant (as opposed to time data).

This paper is organized in a top-down fashion. Definitions and main results are given in Section 2 and in Section 3 we state main upper and lower bounds for the complexity dimensions. After that we concentrate on proving the results; central techniques are discussed in Section 4 and proofs are in Sections 5 and 6. An example of a class with VC-dimension $k$ is given in Section 7.

## 2 Definitions and Statements of Main Results

In this section, we formulate the problem of system identification as a learning problem and illustrate the corresponding learning setting. We define the systems to be studied and, as main results, we state bounds for the number of samples needed in order to identify systems in this setting.

We begin with a general introduction to classification problems. Assume that a set $X$, to be called the *input space*, is given together with a collection $\mathcal{C}$ of mappings $X \to \{0,1\}$.[1] Let $W$ be the set of all sequences

$$w = (u_1, \phi(u_1)), \ldots, (u_s, \phi(u_s))$$

over all $s \geq 1$, $(u_1, \ldots, u_s) \in X^s$ and let $\phi \in \mathcal{C}$. An *identifier* is a map $\psi : W \to \mathcal{C}$. The value of $\psi$ on a sequence $w$ above is denoted as $\psi_w$ instead of $\psi(w)$. The "error" of $\psi$ is the probability that $\psi$ will misclassify a future sample. More formally, the *error of $\psi$ with respect to a probability measure $P$ on $X$, a $\phi \in \mathcal{C}$, and a sequence $(u_1, \ldots, u_s) \in X^s$*, is

$$\text{Err}(P, \phi, u_1, \ldots, u_s) := P\{u \in X \ ; \ \psi_w(u) \neq \phi(u)\}.$$

The class $\mathcal{C}$ is said to be *learnable* if there is some identifier $\psi$ with the following property: For each *accuracy parameter* $\epsilon > 0$ and *confidence parameter* $\delta > 0$ there is some $s$ so that, for every probability $P$ and every $\phi \in \mathcal{C}$,

$$P^s\{(u_1, \ldots, u_s) \in X^s \ ; \ \text{Err}(P, \phi, u_1, \ldots, u_s) > \epsilon\} < \delta,$$

where $P^s$ is the $s$-fold product of $P$. In the learnable case, the function $s(\epsilon, \delta)$ which provides the smallest $s$ achieving this bound for any positive $\epsilon$ and $\delta$ is called the *sample complexity*. It can be proved that learnability is equivalent to the finiteness of the *Vapnik-Chervonenkis dimension* $\text{VC}(\mathcal{C})$ of the class $\mathcal{C}$, which is a combinatorial quantity describing the richness of the class $\mathcal{C}$. Moreover, for learning algorithms that classify the observed samples correctly, the sample complexity is bounded by

$$s(\epsilon, \delta) \leq \max\left\{\frac{8\text{VC}(\mathcal{C})}{\epsilon} \log_2\left(\frac{8e}{\epsilon}\right), \frac{4}{\epsilon} \log_2\left(\frac{2}{\delta}\right)\right\}.$$

---

[1] The set $X$ is assumed to be either countable or an Euclidean space, and the maps in $\mathcal{C}$ are assumed to be measurable. In addition, a mild regularity assumption called "permissibility" is needed so that all sets appearing below are measurable, for further discussion on the topic see an appendix in [15]. In our context the measurability assumptions are satisfied.

In addition, there is a similar lower bound for the sample complexity.

Classification may be viewed as the problem of identifying systems with binary outputs. More generally, we introduce a problem of identification for systems having bounded outputs ($[0,1]$-valued, for technical reasons) via an $L^1$-error, following [5] (for similar statements with $L^2$-error see [2]). Denote by $\mathcal{F}$ a class of mappings from $X$ to $[0,1]$.

By definition, an *identifier* is a mapping from $\cup_{s\in\mathbb{N}} (X \times [0,1])^s$ to $[0,1]^X$. Such a map takes as data a sequence of labeled samples and produces an hypothesis. If $h$ is a $[0,1]$-valued function defined on $X$ and $P$ is a probability measure over $X \times [0,1]$, we define the *error of $h$ with respect to $P$* as

$$\operatorname{Er}_P(h) := \int_{X \times [0,1]} |h(u) - y|\, dP(u,y).$$

For $\epsilon > 0$ and $\delta > 0$ we say that an identifier $\psi$ $(\epsilon, \delta)$-*learns in the agnostic sense with respect to $\mathcal{F}$ from $s$ examples* if, for all distributions $P$ on $X \times [0,1]$, and all $f$ in $\mathcal{F}$:

$$P^s\{w\ ;\ \operatorname{Er}_P(\psi_w) \geq \inf_{f \in \mathcal{F}} \operatorname{Er}_P(f) + \epsilon\} < \delta.$$

Similarly, for $\epsilon > 0$ the function class $\mathcal{F}$ is said to be $\epsilon$-*agnostically learnable* if there is a function $s_0 : (0,1) \to \mathbb{N}$ such that, for all $0 < \delta < 1$, there is an identifier $\psi$ which $(\epsilon, \delta)$-learns in the above sense with $s_0$ samples. In addition, if the identifier always chooses a hypothesis from $\mathcal{F}$, we say that $\mathcal{F}$ is *properly $\epsilon$-agnostically learnable*.

For learning $[0,1]$-valued functions, a sample complexity result may be stated in terms of a fat-shattering dimension, which is a generalization of the VC-dimension. For $\epsilon > 0$ and $\delta > 0$ there is an identifier $\psi$ that properly $(\epsilon, \delta)$-learns in the agnostic sense with respect to $\mathcal{F}$ from ([5])

$$\frac{4}{\alpha^2}\left(\frac{6d}{\ln 2} \ln \frac{7}{\alpha}\left(\frac{336e}{\alpha^3 \ln 2} \ln \frac{7}{\alpha}\right) + \ln \frac{8}{\delta}\right)$$
$$= O\left(\frac{1}{\alpha^2}\left(d \log^2 \frac{1}{\alpha} + \log \frac{1}{\delta}\right)\right)$$

samples, where $0 < \alpha < \epsilon/4$ is chosen so that $d = \operatorname{fat}_{\epsilon/4-\alpha}(\mathcal{F})$ is finite. The quantity $\operatorname{fat}_\gamma(\mathcal{F})$ is called the *fat-shattering dimension* of the class $\mathcal{F}$ and it measures the richness of the class $\mathcal{F}$ with scale $\gamma$.

The sample complexity results show us that the difficulty of system identification in the learning theoretic setting can be analyzed by studying various complexity dimensions, and it is the main focus of this paper. However, formal definitions of the complexity dimensions are delayed until Section 3.

## 2.1 Linear Systems

In the context of learning we discuss continuous-time linear control systems:

$$\dot{x} = Ax + Bu, \qquad x(0) = x^0, \qquad y = Cx, \tag{2.1}$$

where $A$, $B$, and $C$ are $n \times n$, $n \times m$, and $p \times n$ real matrices, and the time interval is [0,1]. We study sign-observations (see [13] for related work in control theory):

$$\text{sign } y(1) = (\text{sign } y_1(1), \ldots, \text{sign } y_p(1))^T,$$

where sign $z = 0$, if $z \leq 0$, sign $z = 1$, if $z > 0$ and $^T$ stands for the transpose. For scalar observations this is a classification problem; each output is classified either 0 or 1 and the VC-dimension can be used to study the complexity of the problem. When $p > 1$, a generalization of the VC-dimension or a loss function is needed.

In general, unlike the VC-dimension associated to discrete-time linear systems [7, 12], the VC-dimension of the classification problem for continuous-time control systems is unbounded [23], even when $n = 1$, and the identification problem is not learnable in the sense discussed earlier. Therefore we restrict the class of admissible controls in order to achieve a bound for the VC-dimension. We consider controls $u = (u_1, \ldots, u_m)$ such that

$$u = G\omega,$$

where $G$ is a $m \times k$ matrix that parameterizes the control. The set of basis input functions $\Omega = \{\omega_1, \ldots, \omega_k\}$ is fixed. The bounds for the VC-dimension or other complexity dimensions will depend on the properties of the set $\Omega$.

For scalar inputs (i.e., $m = 1$) the VC-dimension associated to the mapping from inputs $G$ to scalar sign-observations is bounded by $k$, which in fact can be very large in applications. This bound is tight; we give an example of a function class $\Omega$ for which the associated VC-dimension is indeed $k$ (see Section 7). However, by considering band-limited controls a better bound can be achieved. In this work we consider the following set of basis input functions

$$\Omega = \Big\{\omega_1, \ldots, \omega_k \ ; \ \omega_1, \ldots, \omega_k \text{ linearly independent and}$$
$$\omega_j = t^{\ell_j} e^{\alpha_j t} \sin(\beta_j t) \text{ or } \omega_j = t^{\ell_j} e^{\alpha_j t} \cos(\beta_j t) \quad (2.2)$$
$$\text{with } \ell_j \in \mathbb{N}, \alpha_j, \beta_j \in \mathbb{R}, j = 1, \ldots, k\Big\},$$

and let

$$\ell_{\max} = \max\{\ell_1, \ldots, \ell_k\}. \quad (2.3)$$

The results in this paper hold with minor modifications if basis input functions $\omega_j, j = 1, \ldots, k$ are, for example, linear combinations of the functions of the above form. However, the proofs and formulae are messier.

**Definition 1. (Sign system concept class, $\mathcal{C}_{m,p}$)** Order the set of basis input functions $\Omega$ and denote $\omega = (\omega_1, \ldots, \omega_k)^T$. Let

$$X_\Omega = \{G\omega : [0,1] \to \mathbb{R}^m \ ; \ G \in \mathbb{R}^{mk}\},$$

and for each linear system $\Sigma = (A, B, C, x^0)$ of dimension $n$ define the mapping $\Phi_\Sigma : X_\Omega \to \mathbb{R}^p$ by $\Phi_\Sigma(G\omega) = y(1)$, where $y(1)$ is the solution of $\Sigma$ with control $u = G\omega$. Similarly we define the mapping for sign-observations,

$$S_\Sigma : X_\Omega \to \{0,1\}^p \quad G\omega \mapsto \text{sign}(\Phi_\Sigma(G\omega)).$$

The class of above mappings is the *sign system concept class*, $\mathcal{C}_{m,p} = \{S_\Sigma \; ; \; \Sigma \text{ linear system of dimension } n\}$.

When studying the learning complexity associated to the classification problem we consider the time interval $[0, 1]$ for simplicity, as the length of the interval plays no role in the result. However, for learning $[0, 1]$-valued functions utilizing pseudo-dimension without a loss function, we consider the time interval $[0, \tau]$ with $\tau > 1$.

For the system
$$\dot{x} = Ax + Bu, \qquad (2.4)$$
$$y = Cx,$$

we can write $Ce^{At}B = [\gamma_1, \ldots, \gamma_m]$, where each $\gamma_i$ is a linear combination of $n$ functions $\xi_1, \ldots, \xi_n$. Each $\xi_i$ is of the form $t^\ell e^{at} \sin(bt)$ or $t^\ell e^{at} \cos(bt)$ with $\ell \in \{0, \ldots, n-1\}$ and $a + ib$ an eigenvalue of $A$. Assume that $A$ has a fixed Jordan block structure and let $a_k + ib_k$ be an eigenvalue of $A$. We take $\alpha_{11}, \ldots, \alpha_{nm}, a_1, b_1, \ldots, a_r, b_r$ to be the system parameters, where $\gamma_i(t) = \sum_{j=1}^n \alpha_{ij} \xi_j(t)$ for $i = 1, \ldots, m$ and $a_1, b_1, \ldots, a_r, b_r$ are $n$ eigenvalue parameters. For example the eigenvalue parameters for a real $4 \times 4$ matrix $A$ with eigenvalues $a_1 \pm b_1 i$, $a_2$ and $a_3$ would be $a_1, b_1, a_2$ and $a_3$. Similarly, the eigenvalue parameters with purely complex eigenvalues $a_1 \pm b_1 i$ and $a_2 \pm b_2 i$ would be $a_1, b_1, a_2$ and $b_2$, whereas real eigenvalue parameters would be listed as $a_1, a_2, a_3$ and $a_4$.

Let $U \subset \mathbb{R}^{mk}$ be a bounded set. Define a mapping $F : \lambda \times U \to \mathbb{R}$ such that $F(\lambda, G) = y(\tau)$, where $\tau \geq 1$ and $y(\tau)$ is a solution of (2.4) with system parameters $\lambda$, and initial condition $x(0) = 0$.

In the following definition we take the final time to be $\tau > 1$ in order to show the effect of the time interval in the learning complexity.

**Definition 2.** Assume that the system $\dot{x} = Ax + Bu$, $y = Cx$, $x(0) = 0$ can be parameterized by $\lambda \in \mathbb{R}^{n(m+1)}$ as above and $\|\lambda\|_\infty = \max_{1 \leq i \leq n(m+1)} \lambda_i < 1$. Let $F(\lambda, u) = y(\tau)$ be the solution of (2.4) with system parameters $\lambda$ and control $u = (u_1, \ldots, u_m) \in U = \{u = (u_1, \ldots, u_m) \; ; \; \int_0^\tau u_i(t)dt \leq M, \; i = 1, \ldots, m\}$. Denote $B_\infty^k(c) := \{x \in \mathbb{R}^k \; ; \; \|x\|_\infty < c\}$ and define
$$\mathcal{F}_B = \{F(\lambda, \cdot) : U \to \mathbb{R} \; ; \; \lambda \in B_\infty^{n(m+1)}(1)\}.$$

## 2.2 Main Results

We formulate two theorems about bounding sample complexities as main results. In Section 3 we summarize upper and lower bounds for complexity dimensions studied in this paper together with definitions. The following results, formulated as sample complexity statements, are immediate corollaries of learning complexity bounds proved in this paper, and hence no separate proofs are given.

**Theorem 2.1 (Sample complexity for concept learning).** *For sign systems concept class $\mathcal{C}_{m,1}$ with scalar observations, i.e., $p = 1$, the sample complexity $s(\epsilon, \delta)$ for identifiers that agree with the observed sample can be bounded as*
$$s(\epsilon, \delta) \leq \max \left\{ \frac{8 VC(\mathcal{C}_{m,1})}{\epsilon} \log_2 \left( \frac{8e}{\epsilon} \right), \frac{4}{\epsilon} \log_2 \left( \frac{2}{\delta} \right) \right\},$$

*where*

$$VC(\mathcal{C}_{m,1}) \leq 2(2mn^2 + 4n + 1)\log_2\left[8e\big(8mn^2k(n+\ell_{\max}) + 1\big)\big(2nk + 2(1+2k)^n\big)\right]$$

*and $\ell_{\max}$ is given by (2.3).*

In terms of $n$ (the dimension of the state space) and $k$ (the band-width) the upper bound for the VC-dimension is of the form $O(n^3 \log_2(nk))$. The next section states also a corresponding VC-dimension lower bound, in terms of the band-width, of the form $O(\log(k))$, and together with a lower bound for the sample complexity, this provides an estimate for the number of samples needed in learning. In particular, in a typical setting of fairly small system dimension $n$ and large band-width $k$, the $\log k$ bound is a clear improvement over the linear bound given by elementary analysis.

**Theorem 2.2 (Sample complexity for proper agnostic learning).** *Let $\kappa > 0$, then the class $\mathcal{F}_B$, given in Definition 2, is properly agnostically learnable from*

$$O\left(\frac{1}{\epsilon^2}\left(\mathrm{fat}_{(1/4-\kappa)\epsilon}(\mathcal{F}_B)\log^2\frac{1}{\epsilon} + \log\frac{1}{\delta}\right)\right)$$

*samples, where*

$$\mathrm{fat}_{(1/4-\kappa)\epsilon}(\mathcal{F}_B) \leq \min\begin{cases}(m+1)n\log_2\left\lfloor\frac{n^2 m\tau^n e^\tau kM}{(1/4-\kappa)\epsilon}\right\rfloor, \\ 2(m+4)n\log_2\left(8e\big(nmk4(n+\ell_{\max})+1\big)\big(2nk+2(2k+1)^n\big)\right),\end{cases}$$

*together with $\ell_{\max}$ given by (2.2) and (2.3), and $M$ a constant satisfying*

$$\int_0^\tau |u_i(\tau - t)|dt \leq kM$$

*for all $i = 1, \ldots, m$. In above, $\lfloor x \rfloor$ stands for the integer part of $x$.*

In addition to the bound above the control systems can be parameterized linearly as follows: Let $u(s) = \sum_{j=1}^k g_j \omega_j(s)$. Then

$$y(\tau) = C\int_0^\tau e^{A(\tau-s)}Bu(s)ds = C\int_0^\tau e^{A(\tau-s)}B\sum_{j=1}^k g_j\omega_j(s)ds$$
$$= \sum_{j=1}^k g_j \underbrace{C\int_0^\tau e^{A(\tau-s)}B\omega_j(s)ds}_{\lambda_j(\tau)} = \sum_{j=1}^k g_j\lambda_j(\tau). \quad (2.5)$$

Then by considering the class $\mathcal{F} := \{g \in \mathbb{R}^k \mapsto \sum_{j=1}^k g_j\lambda_j(\tau)\,;\, \sum_{j=1}^k \lambda_j(\tau)^2 = 1\}$ and further restricting $g$ so that $\|g\|_2 \leq R$, [18, 16] have shown that

$$\mathrm{fat}_\gamma(\mathcal{F}) \leq \min\{9R^2/\gamma^2, k+1\} + 1.$$

# 3 Complexity dimensions; main upper and lower bounds

We begin this section by defining various learning complexity dimensions and after that we summarize the main upper and lower bounds proved in this paper.

**Definition 3 (Vapnik-Chervonenkis dimension).** The richness of the collection $\mathcal{C}$ can be measured by its Vapnik-Chervonenkis (VC) dimension introduced in [19]. A set $S = \{x_1, \ldots, x_n\} \subseteq X$ is said to be *shattered* by $\mathcal{C}$ if, for every subset $B \subseteq S$, there exists a set $A \in \mathcal{C}$ such that $S \cap A = B$. The *Vapnik-Chervonenkis dimension of $\mathcal{C}$*, denoted $\mathrm{VC}(\mathcal{C})$, equals the largest integer $n$ such that there exists a set of cardinality $n$ that is shattered by $\mathcal{C}$. For example, in $\mathbb{R}^k$ the VC-dimension of closed half-spaces through the origin is $k$ [22]. Thus, if $\mathrm{VC}(\mathcal{C}) = d$, $\mathcal{C}$ is not rich enough to distinguish all subsets of any $d+1$ element set, but there is some $d$ element set where subsets can be distinguished. Proving exact values of the VC-dimension is hard and typically one looks for upper and lower bounds for the VC-dimension, as is also done in this paper.

For our purposes, it is more convenient to work with shattering in terms of dichotomies, i.e., boolean-valued maps. We identify subsets of $D$ with Boolean functions $\phi : D \to \{0,1\}$. Similarly, each set $C \in \mathcal{C}$ gives rise to a Boolean function on $X$, and intersections $C \cap D$ are restrictions of functions to $D$. In this language, a subset $D \subset X$ is shattered by $\mathcal{F} := \{\phi \; ; \; \phi : X \to \{0,1\}\}$, if every dichotomy on $D$ is a restriction to $D$ of some $\phi \in \mathcal{F}$.

**Definition 4 (Pseudo-dimension with respect to a loss function).** For a given class of functions $\mathcal{F} : X \to Y$ and a loss function $L : Y \times Y \to [0, r]$, we introduce for each $f \in \mathcal{F}$ the function

$$A_{f,L} : X \times Y \times \mathbb{R} \to \{0,1\}; \quad (x, y, \rho) \mapsto \mathrm{sign}(L(f(x), y) - \rho), \tag{3.1}$$

and let $A_{\mathcal{F},L}$ denote all such $A_{f,L}$ with $f \in \mathcal{F}$. The *pseudo-dimension of $\mathcal{F}$ with respect to the loss function $L$*, $\mathrm{PD}[\mathcal{F}, L]$, is defined as:

$$\mathrm{PD}[\mathcal{F}, L] := \mathrm{VC}(A_{\mathcal{F},L}).$$

The VC-dimension characterizes learnability of $\{0,1\}$-valued functions. For learning real-valued functions we look for a generalization of the VC-dimension with similar properties. One such generalization is the pseudo-dimension. Unfortunately, pseudo-dimension does not share the property the VC-dimension has; there are learnable function classes with infinite pseudo-dimension, see [20, p.206] and [3].

Next we define the fat-shattering dimension that corresponds to shattering with fixed "margin" $\gamma$. Both the pseudo-dimension and the fat-shattering dimension can be used to bound certain covering numbers and in this sense they act like the VC-dimension. Moreover, the fat-shattering dimension gives upper and lower bounds for covering numbers of function classes and the finiteness of the fat-shattering dimension can characterize learnability (see [1] and [2]).

**Definition 5 (Fat-shattering dimension).** Let $\mathcal{F}$ be a set of real-valued functions. We say that a set of points $X$ is $\gamma$-*shattered by* $\mathcal{F}$ if there are real numbers $r_x$ indexed by $x \in X$ such that for all binary vectors $b_x$ indexed by $X$, there is a function $f_b \in \mathcal{F}$ satisfying

$$f_b(x) \geq r_x + \gamma, \quad \text{if } b_x = 1, \text{ and}$$
$$f_b(x) \leq r_x - \gamma, \quad \text{otherwise.}$$

The *fat-shattering dimension* $\text{fat}_\gamma(\mathcal{F})$ is a function from positive real numbers to integers which maps a value $\gamma$ to the size of the largest fat-shattered set, if it is finite, or infinity otherwise.

The shattering dimension when the margin $\gamma$ equals 0 is called the *pseudo-dimension* and it is denoted by $\text{PD}(\mathcal{F})$. Clearly, for all $\gamma > 0$, $\text{fat}_\gamma(\mathcal{F}) \leq \text{PD}(\mathcal{F})$.

## 3.1 Bounds

We begin by stating bounds in the easiest learning setting—classifying the final state observations as either 0 or 1.

**Theorem 3.1. (VC-dimension upper bound, p = 1)** *The VC-dimension of the sign system concept class $\mathcal{C}_{m,1}$ with scalar observations can be bounded as*

$$VC(\mathcal{C}_{m,1}) \leq 2(2mn^2 + 4n + 1)\log_2\left[8e\big(8mn^2k(n + \ell_{\max}) + 1\big)\big(2nk + 2(1 + 2k)^n\big)\right],$$

*where $\ell_{\max}$ is given by* (2.3).

In terms of $n$ (the dimension of the state space) and $k$ (the band-width) the upper bound is of the form $O(n^3 \log_2(nk))$.

All VC-dimension upper bounds are based on the fact that input basis functions satisfy a certain rationality condition. Remark 5.3 indicates how the bound is formed when the input functions satisfy the more abstract rationality condition. In that case the degrees of the polynomials and the number of polynomial evaluations are different. However, in terms of $n$ and $k$, the bound is of the same form. VC-dimension or pseudo-dimension bounds stated in this paper can be modified for the rationality condition in the same way.

The lower bound for the VC-dimension is in terms of $n$ and $k$. It holds for linearly independent continuous basis input functions and compared to upper bounds, no particular form of the functions is needed. The bound is obtained by imposing a specific structure on control systems, and a lower bound for a restricted class of control systems provides a lower bound for more general classes.

**Theorem 3.2. (VC-dimension lower bound, m = 1, p = 1)**

$$VC(\mathcal{C}_{1,1}) \geq \max\left\{m'\left\lfloor\log_2\left\lfloor\frac{k}{m'}\right\rfloor\right\rfloor, m'\right\}$$

*where $m' = \min\{n, k\}$.*

In terms of $k$ the upper and the lower bound match up to a constant. For $n$ and $k$ the lower bound is typically of the form $O(n \log_2(k/n))$. Note that if the system dimension $n$ is small compared to the band-width $k$ the VC-dimension upper and lower bounds in Theorems 3.1 and 3.2 become tighter, both being of the form $c \log_2 k$ (with different values of the constant $c$).

Extending the upper bounds to the case of vector-valued observations can be done in various ways based on the result obtained for scalar observations. For example, we may consider the $p$-dimensional output as bits representing a number in $\{0, \ldots, 2^p - 1\}$ and introduce a loss function for each $f \in \mathcal{C}_{m,p}$ as $L_{0\text{-}1,f}(z, a) = L_{0\text{-}1}(f(z), a) = 1$, when $f(z) \neq a$, and 0 otherwise. We define the VC-dimension of the $p$-dimensional observation as the VC-dimension of the above class of loss functions. Modifying the argument used with scalar observations leads to a bound of the following form:

**Theorem 3.3. (VC-dimension upper bound)**

$$VC(\mathcal{C}_{m,p}) \leq 2(2pmn^2 + 4n + p)$$
$$\times \log_2 \Big[ 8e \big( 8mn^2 k(n + \ell_{\max}) + 1 \big) \big( 2^p - 1 + 2p(2k+1)^n + 2nk \big) \Big],$$

where $\ell_{\max}$ is given by (2.3).

Next we state the main result concerning learnability of the actual input-output mapping, i.e., learning without taking the sign of the final state observation.

**Definition 6. (Control system concept class, $\mathcal{G}_{p,L}$)** Let $\widetilde{\mathcal{F}} = \{\Phi_\Sigma : X_\Omega \to \mathbb{R}^p \; ; \; \Sigma \text{ linear system of dimension } n\}$ and define the *control system concept class* as $\mathcal{G}_{p,L} = A_{\widetilde{\mathcal{F}},L}$, where $A_{\widetilde{\mathcal{F}},L}$ is given by (3.1).

Methods for calculating upper bounds for the VC-dimension readily give a tool for obtaining upper bounds for the pseudo-dimension with respect to loss functions that preserve the rationality structure of the output. A typical example is illustrated by the loss function

$$L(z_1, z_2) = (z_1 - z_2)^2 / (1 + (z_1 - z_2)^2) \tag{3.2}$$

and the following result:

**Theorem 3.4. (Pseudo-dimension upper bound, p = 1)**

$$PD(\mathcal{G}_{1,L}) \leq 2 \left( 2mn^2 + 4n + 1 \right) \log_2 \Big[ 16e \big( 8mn^2 k(n + \ell_{\max}) + 1 \big) \big( 2nk + 2(2k+1)^n \big) \Big],$$

where the loss function $L$ is given by (3.2) and $\ell_{\max}$ is given by (2.3).

This differs from the corresponding VC-dimension bound only by the maximum degree of the polynomials, which is doubled. Extending this pseudo-dimension bound for $p$-dimensional observations can be done naturally by modifying the loss function. Lower bounds for the VC-dimension are lower bounds for the pseudo-dimension as such.

The next results summarize upper bounds for the fat-shattering dimension. We begin by illustrating how fat-shattering dimension can be bounded for Lipschitz functions in certain cases:

**Theorem 3.5 (Fat-shattering bound).** *Let $F(\lambda, u) : \mathbb{R}^k \times U \to \mathbb{R}$ be such that $F(\cdot, u)$ is Lipschitz with constant $L$, i.e., $|F(\lambda_1, u) - F(\lambda_2, u)| \leq L\|\lambda_1 - \lambda_2\|$ for all $u \in U$. For any subset $B \subseteq \mathbb{R}^k$, consider the following class of functions:*

$$\mathcal{F}_B = \{F(\lambda, \cdot) : U \to \mathbb{R} \; ; \; \lambda \in B\}.$$

*Then*

$$\mathrm{fat}_\gamma(\mathcal{F}_{B_\infty^k(C)}) \leq k \log_2 \left\lfloor \frac{CL}{\gamma} \right\rfloor,$$

$$\mathrm{fat}_\gamma(\mathcal{F}_{\bar{B}_\infty^k(C)}) \leq k \log_2 \left(1 + \left\lfloor \frac{CL}{\gamma} \right\rfloor\right),$$

$$\mathrm{fat}_\gamma(\mathcal{F}_{\bar{B}_{2,1}}) \leq k \log_2 \left(C + \frac{L}{\gamma}\right),$$

*where*

$$B_\infty^k(C) = \{x \in \mathbb{R}^k \; ; \; \|x\|_\infty = \max_{1 \leq i \leq k} |x_i| < C\},$$

$$\bar{B}_\infty^k(C) = \{x \in \mathbb{R}^k \; ; \; \|x\|_\infty = \max_{1 \leq i \leq k} |x_i| \leq C\},$$

$$\bar{B}_2^k(C) = \{x \in \mathbb{R}^k \; ; \; \|x\|_2 = \sqrt{\sum_{i=1}^k x_i^2} \leq C\}.$$

**Theorem 3.6 (Fat-shattering bound for a control system).** *Assume that the system $\dot{x} = Ax + Bu$, $y = Cx$, $x(0) = 0$ can be parameterized by $\lambda \in \mathbb{R}^{n(m+1)}$ as in Definition 2 and $\|\lambda\|_\infty < 1$. Let $F(\lambda, u) = y(\tau)$ be the solution with system parameters $\lambda$ and control $u = (u_1, \ldots, u_m) \in U = \{u = (u_1, \ldots, u_m) \; ; \; \int_0^\tau u_i(t)dt \leq M, i = 1, \ldots, m\}$. Define*

$$\mathcal{F}_B = \{F(\lambda, \cdot) : U \to \mathbb{R} \; ; \; \lambda \in B_\infty^k(1)\}.$$

*Then*

$$\mathrm{fat}_\gamma(\mathcal{F}_B) \leq n(m+1) \log_2 \left\lfloor \frac{n^2 m \tau^n e^\tau M}{\gamma} \right\rfloor.$$

# 4 Techniques for Proving VC-Dimension Results

Our main results are based on a fact that the basis input functions satisfy a certain rationality condition. In this section we first formulate this rationality condition and then we summarize existing results that are used in proving upper and lower bounds for the complexity dimensions.

We recall briefly the control system setting. We study systems

$$\dot{x} = Ax + Bu, \quad x(0) = x^0, \quad y = Cx,$$
$$u = G\omega, \quad \text{and} \quad \omega_j \in \Omega \text{ for } j = 1, \ldots, k,$$

with basis input functions

$$\Omega = \Big\{\omega_1, \ldots, \omega_k \, ; \omega_1, \ldots, \omega_k \text{ linearly independent and}$$
$$\omega_j = t^{\ell_j} e^{\alpha_j t} \sin(\beta_j t) \text{ or } \omega_j = t^{\ell_j} e^{\alpha_j t} \cos(\beta_j t)$$
$$\text{with } \ell_j \in \mathbb{N}, \alpha_j, \beta_j \in \mathbb{R}, j = 1, \ldots, k\Big\},$$

such that

$$\ell_{\max} = \max\{\ell_1, \ldots, \ell_k\}.$$

**Definition 7. (Rationality condition (RAT))** Let $n$ be a positive integer. We say that a bounded function $\omega : [0, 1] \to \mathbb{R}$ satisfies the rationality condition relative to the class of $n$-dimensional systems if there exists $h$ polynomial functions $f_1, \ldots, f_h : \mathbb{R}^4 \to \mathbb{R}$ and $2\gamma n$ rational functions $r_{ij\tilde{\ell}}$, $i \in \{1, 2\}$, $j \in \{1, \ldots, \gamma\}$ and $\tilde{\ell} \in \{1, \ldots, n\}$ with no poles on subsets $S_{ij\tilde{\ell}}$ of $\mathbb{R}^4$, such that the following properties hold:

1. For each $i \in \{1, 2\}$, $\tilde{\ell} \in \{1, \ldots, n\}$, $\mathbb{R}^4$ is a disjoint union of $S_{i1\tilde{\ell}}, \ldots, S_{i\gamma\tilde{\ell}}$.

2. Each $S_{ij\tilde{\ell}}$ can be defined in terms of a Boolean expression involving $[f_1 = 0], \ldots, [f_h = 0]$, where we say that for functions $f_1, \ldots, f_h : \mathbb{R}^4 \to \mathbb{R}$, $[f_i = 0]$ has value 1, if $f_i(x_1, x_2, x_3, x_4) = 0$, and 0 otherwise.

3. Letting $r_{i\tilde{\ell}} : \mathbb{R}^4 \to \mathbb{R}$, $i \in \{1, 2\}$, $\tilde{\ell} \in \{1, \ldots, n\}$ be defined as

$$r_{i\tilde{\ell}}(v) = \begin{cases} r_{i1\tilde{\ell}}(v), & \text{if } v \in S_{i1\tilde{\ell}}, \\ \vdots & \vdots \\ r_{i\gamma\tilde{\ell}}(v), & \text{if } v \in S_{i\gamma\tilde{\ell}}, \end{cases}$$

then for each $a, b \in \mathbb{R}$, and for all $\tilde{\ell} \in \{1, \ldots, n\}$

$$\int_0^1 t^{\tilde{\ell}-1} e^{at} \cos(bt)\omega(t)dt = r_{1\tilde{\ell}}(a, b, e^a \cos b, e^a \sin b),$$
$$\int_0^1 t^{\tilde{\ell}-1} e^{at} \sin(bt)\omega(t)dt = r_{2\tilde{\ell}}(a, b, e^a \cos b, e^a \sin b).$$

We denote by $d_{\max}$ the maximum degree of any polynomial (i.e., $f_1, \ldots, f_h$, numerators and denominators of $r_{ij\tilde{\ell}}$'s) appearing in the rationality condition .

*Remark 4.1.* First, entries of $e^{At}$ are functions of the form $t^s e^{at} \cos(bt)$ and $t^s e^{at} \sin(bt)$. Solving (2.1) involves convolutions of $e^{At}$ and the basis input functions $\omega_j$, and we require those to be rational functions.

*Example 4.2.* Let $\omega(t) = \sin(ct)$, with nonzero $c$. Then

$$\int_0^1 e^{at} \sin(bt)\omega(t)dt = \frac{1}{2}\int_0^1 e^{at}\cos((b-c)t)dt - \frac{1}{2}\int_0^1 e^{at}\cos((b+c)t)dt$$

after integration this can be split into cases with no poles yielding

$$\int_0^1 e^{at}\sin(bt)\omega(t)dt = \begin{cases} \frac{p(a,b,e^a\cos b, e^a \sin b)}{(a^2+(b+c)^2)(a^2+(b-c)^2)}, & \text{if } f_1 \neq 0, f_2 \neq 0, \\ \frac{b-\sin b \cos b}{2b}, & \text{if } f_1 = 0, f_2 \neq 0, \\ \frac{\sin b \cos b - b}{2b}, & \text{if } f_1 \neq 0, f_2 = 0, \end{cases}$$

where

$$f_1(a,b,e^a\sin b, e^a\cos b) = a^2 + (b+c)^2,$$
$$f_2(a,b,e^a\sin b, e^a\cos b) = a^2 + (b-c)^2,$$

and $p(a,b,e^a\cos b, e^a\sin b)$ stands for the polynomial

$$\begin{aligned} &-4abc + 4abce^a\cos b\cos c - 2a^2ce^a\cos c\sin b + 2b^2ce^a\cos c\sin b \\ &- 2c^3e^a\cos c\sin b - 2a^2be^a\cos b\sin c - 2b^3e^a\cos b\sin c + 2bc^2e^a\cos b\sin c \\ &+ 2a^3e^a\sin b\sin c + 2ab^2e^a\sin b\sin c + 2ac^2e^a\sin b\sin c. \end{aligned}$$

**Lemma 4.3.** *Each $\omega_j \in \Omega$ given by (2.2) satisfies the rationality condition. Further, the maximum degree of polynomials in (RAT) is at most $4(n+\ell_{\max})$, where $\ell_{\max}$ is given by (2.3).*

### Review of VC-Dimension Techniques

In the context of control theory it is sometimes easier to work with the dual VC-dimension. Assume that a function $F : \Lambda \times X \to \{0,1\}$ is given. This induces two function classes

$$\mathcal{F} := \{F(\lambda, \cdot) : X \to \{0,1\} \; ; \; \lambda \in \Lambda\}$$

and

$$\mathcal{F}^* := \{F(\cdot, x) : \Lambda \to \{0,1\} \; ; \; x \in X\}.$$

The complexity dimension $\text{VC}(\mathcal{F}^*)$ is called the *dual VC-dimension* of $\mathcal{F}$ and it is related to $\text{VC}(\mathcal{F})$ as follows [20]:

$$\text{VC}(\mathcal{F}) \geq \lfloor \log_2 \text{VC}(\mathcal{F}^*) \rfloor, \tag{4.1}$$

where $\lfloor x \rfloor$ is the integer part of $x$.

A sharper estimate can be obtained if $\Lambda$ can be written as a product $\Lambda_1 \times \cdots \times \Lambda_n$. The following construction and result are due to DasGupta and Sontag [7]. We study in particular those dichotomies that are defined on "rectangular" subsets of $\Lambda$. Let $L = L_1 \times \cdots \times L_n$ be a subset of $\Lambda$ such that for each $i$, $L_i \subset \Lambda_i$ is nonempty. Given any index $1 \leq \kappa \leq n$ a *$\kappa$-axis dichotomy* on $L$ is any function $\delta : L \to \{0, 1\}$ which depends only on the $\kappa$th coordinate, i.e., there is some function $\phi : L_\kappa \to \{0, 1\}$ so that $\delta(\lambda_1, \ldots, \lambda_n) = \phi(\lambda_\kappa)$ for all $(\lambda_1, \ldots, \lambda_n) \in L$. We say that a mapping is an axis dichotomy if it is a $\kappa$-axis dichotomy for some $\kappa$. A rectangular set $L$ is said to be *axis shattered* by $\mathcal{F}^*$ if every axis dichotomy is a restriction to $L$ of some function of the form $F(\cdot, x) : \Lambda \to \{0, 1\}$, for some $x \in X$.

**Theorem 4.4. (Axis shattering bound [7])** *If $L = L_1 \times \cdots \times L_n \subset \Lambda$ can be axis shattered and each $L_i$ has cardinality $r_i > 0$, then $VC(\mathcal{F}) \geq \lfloor \log_2(r_1) \rfloor + \cdots + \lfloor \log_2(r_n) \rfloor$.*

Upper bounds for VC-dimensions of concept classes that are obtained by evaluating polynomial equalities and inequalities can be obtained in terms of the number and degrees of the polynomials:

**Theorem 4.5. (Goldberg-Jerrum bound [10])** *Given a function $F : \Lambda \times X \to \{0, 1\}$ and the associated concept class $\mathcal{F} := \{F(\lambda, \cdot) : X \to \{0, 1\} \ ; \ \lambda \in \Lambda\}$, suppose that $\Lambda = \mathbb{R}^\ell$ and $X = \mathbb{R}^k$. Let $F$ be defined in terms of a Boolean formula involving at most $s$ polynomial equalities and inequalities in $\ell + k$ variables, each polynomial being of degree at most $d$ in $\lambda$ for all $x \in \mathbb{R}^k$. Then, $VC(\mathcal{F}) \leq 2\ell \log_2(8eds)$.*

The Goldberg-Jerrum bound is based on a result showing that the number of sign assignments $\{-1, 0, 1\}$ to polynomials can not grow too fast:

**Theorem 4.6. ([10])** *Suppose that $f_1, \ldots, f_m$ are polynomials of degree at most $d$ in $n \leq m$ variables. Then the number of distinct vectors*

$$\big[\text{sign } f_1(x), \ldots, \text{sign } f_m(x)\big] \in \{-1, 0, 1\}^m,$$

*that can be generated by varying $x$ over $\mathbb{R}^n$, is at most $((8edm)/n)^n$.*

## 5 Proofs of the VC-dimension bounds

### 5.1 An Upper Bound for the VC-Dimension with Scalar Observations

We begin this section by proving Lemma 4.3 stating that the input basis functions satisfy the rationality condition (RAT) and bounding the degrees of polynomials appearing in (RAT). As a proposition we formalize how control systems can be parameterized. After that, as a lemma, we develop an upper bound for the VC-dimension induced by the control system (2.1) with its initial state fixed to be zero. Theorem 3.1 with an arbitrary initial condition is then a simple modification of the argument.

**Proof of Lemma 4.3.** If $\omega(t) = t^\ell e^{\alpha t} \sin(\beta t)$ or $\omega(t) = t^\ell e^{\alpha t} \cos(\beta t)$ with $\ell \leq \ell_{\max}$ then in place of

$$\int_0^1 t^{\tilde{\ell}} e^{at} \sin(bt) \omega(t) dt \quad \text{or} \quad \int_0^1 t^{\tilde{\ell}} e^{at} \cos(bt) \omega(t) dt \tag{5.1}$$

by combining exponents and using sum formulae for sin and cos (see Example 4.2) it is enough to study terms of the form

$$\int_0^1 t^{\tilde{k}} e^{\tilde{a}t} \sin(\tilde{b}t) dt \quad \text{or} \quad \int_0^1 t^{\tilde{k}} e^{\tilde{a}t} \cos(\tilde{b}t) dt,$$

where $\tilde{k} \in \{0, \ldots, n + \ell_{\max} - 1\}$. In fact, each expression in (5.1) is one of the following type

$$\frac{1}{2} \left( \int_0^1 t^{\tilde{\ell}} e^{\tilde{a}t} \cos\left((b - \beta)t\right) dt - \int_0^1 t^{\tilde{\ell}} e^{\tilde{a}t} \cos\left((b + \beta)t\right) dt \right)$$

$$\frac{1}{2} \left( \int_0^1 t^{\tilde{\ell}} e^{\tilde{a}t} \cos\left((b + \beta)t\right) dt + \int_0^1 t^{\tilde{\ell}} e^{\tilde{a}t} \cos\left((b - \beta)t\right) dt \right)$$

$$\frac{1}{2} \left( \int_0^1 t^{\tilde{\ell}} e^{\tilde{a}t} \sin\left((b - \beta)t\right) dt - \int_0^1 t^{\tilde{\ell}} e^{\tilde{a}t} \sin\left((b + \beta)t\right) dt \right),$$

$$\frac{1}{2} \left( \int_0^1 t^{\tilde{\ell}} e^{\tilde{a}t} \sin\left((b + \beta)t\right) dt + \int_0^1 t^{\tilde{\ell}} e^{\tilde{a}t} \sin\left((b - \beta)t\right) dt \right),$$

where $\tilde{a} = a + \alpha$. Because for $\tilde{k} > 0$

$$\int_0^1 t^{\tilde{k}} e^{\tilde{a}t} \sin(\tilde{b}t) dt = \frac{e^{\tilde{a}}}{\tilde{a}^2 + \tilde{b}^2} (\tilde{a} \sin \tilde{b} - \tilde{b} \cos \tilde{b})$$

$$- \frac{\tilde{k}}{\tilde{a}^2 + \tilde{b}^2} \int_0^1 t^{\tilde{k}-1} e^{\tilde{a}t} (\tilde{a} \sin(\tilde{b}t) - \tilde{b} \cos(\tilde{b}t)) dt \quad (5.2)$$

and

$$\int_0^1 e^{\tilde{a}t} \sin(\tilde{b}t) dt = \frac{e^{\tilde{a}}(\tilde{a} \sin \tilde{b} - \tilde{b} \cos \tilde{b}) + \tilde{b}}{\tilde{a}^2 + \tilde{b}^2} \quad (5.3)$$

and a similar formulae for $\int_0^1 t^{\tilde{k}} e^{\tilde{a}t} \cos(\tilde{b}t) dt$ hold, we see by induction that the numerator of $\int_0^1 t^{\tilde{k}} e^{\tilde{a}t} \sin(\tilde{b}t) dt$ is a polynomial of $\tilde{a}, \tilde{b}, e^{\tilde{a}} \cos \tilde{b}$ and $e^{\tilde{a}} \sin \tilde{b}$. By using sum formulae for sin and cos, the previous expression is in turn a polynomial of $a, b, e^a \cos b$ and $e^a \sin b$ because $\tilde{a} = a + \alpha$ and $\tilde{b} = b \pm \beta$ for some fixed $\alpha$ and $\beta$. By similar arguments, the denominator is a polynomial of $a$ and $b$. Note that, for example, $e^\alpha$ equals a constant times $e^a$, so this process does not change the degrees of the polynomials.

Further, observe that the denominator of $\int_0^1 t^{\tilde{\ell}} e^{at} \sin(bt) \omega(t) dt$ consists of at most two products of variables $a$ and $b$ of the form $((a + \alpha)^2 + (b \pm \beta)^2)^{\tilde{\ell}+\ell+1}$, and similarly with the $\cos(bt)$ term. Let us index the basis input functions $\omega_1, \ldots, \omega_k$ so that $\omega_\kappa$ has parameters $\alpha_\kappa$ and $\beta_\kappa$. Hence the functions $f_i$ in (RAT), defining the subsets without poles, can be taken as

$$\{(a + \alpha_\kappa)^2 + (b - \beta_\kappa)^2, (a + \alpha_\kappa)^2 + (b + \beta_\kappa)^2 \; ; \; \kappa = 1, \ldots, k\}.$$

Furthermore, the sets $S_{ij\tilde{\ell}}$ are as simple as

$$\cup_{\kappa=1}^k \{\{(x_1, x_2, x_3, x_4) \; ; \; x_1 = -\alpha_\kappa, x_2 = -\beta_\kappa\} \cup \{(x_1, x_2, x_3, x_4) \; ; \; x_1 = -\alpha_\kappa, x_2 = \beta_\kappa\}\}$$

and

$$\mathbb{R}^4 \setminus \cup_{\kappa=1}^{k} \{\{(x_1, x_2, x_3, x_4) \; ; \; x_1 = -\alpha_\kappa, x_2 = -\beta_\kappa\} \cup \\ \{(x_1, x_2, x_3, x_4) \; ; \; x_1 = -\alpha_\kappa, x_2 = \beta_\kappa\}\}.$$

We turn to estimating the maximum degree of polynomials appearing in (RAT). We already saw that functions $f_i$ are polynomials of degree 2. Equations (5.2) and (5.3) show that the degree of numerator is not higher than the one of denominator. We claim that

$$\int_0^1 t^{\tilde{k}} e^{\tilde{a}t} \left\{ \begin{array}{c} \sin(\tilde{b}t) \\ \cos(\tilde{b}t) \end{array} \right. dt = \frac{P(2(\tilde{k}+1))}{(\tilde{a}^2 + \tilde{b}^2)^{\tilde{k}+1}} \quad \text{for } \tilde{k} = 0, 1, \ldots$$

where $P(2(\tilde{k}+1))$ stands for some polynomial in $\tilde{a}$, $\tilde{b}$, $e^{\tilde{a}} \sin(\tilde{b})$ and $e^{\tilde{a}} \cos(\tilde{b})$ of degree $2(\tilde{k}+1)$. Clearly, the claim is true for $\tilde{k} = 0$ by (5.3) and the inductive argument follows from (5.2). Assuming the claim true for $\tilde{k} - 1$, we get

$$\int_0^1 t^{\tilde{k}} e^{\tilde{a}t} \sin(\tilde{b}t) dt$$

$$= \frac{P(2)}{\tilde{a}^2 + \tilde{b}^2} - \frac{\tilde{k}\tilde{a}}{\tilde{a}^2 + \tilde{b}^2} \int_0^1 t^{\tilde{k}-1} e^{\tilde{a}t} \sin(\tilde{b}t) dt + \frac{\tilde{k}\tilde{b}}{\tilde{a}^2 + \tilde{b}^2} \int_0^1 t^{\tilde{k}-1} e^{\tilde{a}t} \cos(\tilde{b}t) dt$$

$$= \frac{P(2)}{\tilde{a}^2 + \tilde{b}^2} - \frac{P(1)}{(\tilde{a}^2 + \tilde{b}^2)} \frac{P(2\tilde{k})}{(\tilde{a}^2 + \tilde{b}^2)^{\tilde{k}}} + \frac{P(1)}{(\tilde{a}^2 + \tilde{b}^2)} \frac{P(2\tilde{k})}{(\tilde{a}^2 + \tilde{b}^2)^{\tilde{k}}}$$

$$= \frac{P(2)(\tilde{a}^2 + \tilde{b}^2)^{\tilde{k}} - 2P(2\tilde{k}+1)}{(\tilde{a}^2 + \tilde{b}^2)^{\tilde{k}+1}} = \frac{P(2(\tilde{k}+1))}{(\tilde{a}^2 + \tilde{b}^2)^{\tilde{k}+1}}$$

and similarly for the $\cos(\tilde{b}t)$ term, concluding the proof of the claim. As a corollary of the claim

$$\int_0^1 t^{\tilde{\ell}} e^{at} \sin(bt) \omega(t) dt = \frac{P(2(\tilde{k}+1))}{(\tilde{a}^2 + (b+\beta)^2)^{\tilde{k}+1}} + \frac{P(2(\tilde{k}+1))}{(\tilde{a}^2 + (b-\beta)^2)^{\tilde{k}+1}},$$

where $\tilde{k} = \ell + \tilde{\ell}$ and $\tilde{a} = a + \alpha$.

Hence the maximum degree of denominators of expressions in (5.1) is $2(\tilde{k}+1)+2(\tilde{k}+1) = 4(\tilde{k}+1)$ with $\tilde{k} \in \{0, \ldots, \ell_{\max} + n - 1\}$. Thus the maximum degree of polynomials appearing in the (RAT) is $4(n + \ell_{\max})$. $\square$

The next proposition indicates how control systems are parameterized and later the concept or function classes associated to control systems are obtained by varying the parameter vector.

**Proposition 5.1.** *Denote the basis input functions by $\omega = (\omega_1, \ldots, \omega_k)^T$ and assume that each $\omega_i$, $i = 1, \ldots, k$ satisfies the rationality condition (RAT) and let $\Lambda = \mathbb{R}^{2pn^2m} \times \mathbb{R}^{4n} \times \mathbb{R}^p$. Then there exists a mapping $H : \Lambda \times \mathbb{R}^{mk} \to \mathbb{R}^p$ (depending on $\omega$) such that for each $\Sigma = (A, B, C, x^0)$ there exists a $\lambda \in \Lambda$ satisfying*

$$\Phi_\Sigma(G\omega) = H(\lambda, G) \text{ for all } G \in \mathbb{R}^{mk}.$$

*Proof.* Given a system $\Sigma = (A, B, C, x^0)$,

$$\Phi_\Sigma(u) = y(1) = Ce^A x^0 + C \int_0^1 e^{A(1-t)} Bu(t)dt.$$

By an argument based on the real Jordan form of $e^{At}$, the entries of $e^{A(1-t)}$ are linear combinations of functions of the form $t^{\tilde{\ell}} e^{at} \cos(bt)$ and $t^{\tilde{\ell}} e^{at} \sin(bt)$, where $\tilde{\ell} \in \{0, \ldots, n-1\}$ and $a + ib$ is an eigenvalue of $A$. Hence we define the following $2n$ functions

$$\xi_1(a, b, t) = e^{at} \cos(bt)$$
$$\xi_2(a, b, t) = te^{at} \cos(bt)$$
$$\vdots$$
$$\xi_n(a, b, t) = t^{n-1} e^{at} \cos(bt)$$
$$\xi_{n+1}(a, b, t) = e^{at} \sin(bt)$$
$$\vdots$$
$$\xi_{2n}(a, b, t) = t^{n-1} e^{at} \sin(bt).$$

By the rationality condition (RAT), for all $\ell = 1, \ldots, 2n$

$$\int_0^1 \xi_\ell(a, b, t) \omega_j(t) dt = \frac{\hat{P}_{\ell j}(a, b, e^a \cos b, e^a \sin b)}{\hat{Q}_{\ell j}(a, b, e^a \cos b, e^a \sin b)}$$

for all $a, b \in \mathbb{R}$ and where $\hat{P}_{\ell j}$ and $\hat{Q}_{\ell j}$ are piecewise polynomial expressions.

Let $H(\mathbf{A}, \mathbf{X}, h, G) = (H_1, \ldots, H_p)^T$, where for $1 \leq \kappa \leq p$

$$H_\kappa(\mathbf{A}, \mathbf{X}, h, G) = \sum_{i=1}^m \sum_{r=1}^n \sum_{\ell=1}^{2n} \alpha_{ir\ell\kappa} \sum_{j=1}^k g_{ij} \frac{\hat{P}_{\ell j}(x_{r1}, x_{r2}, x_{r3}, x_{r4})}{\hat{Q}_{\ell j}(x_{r1}, x_{r2}, x_{r3}, x_{r4})} + h_\kappa$$

and

$$\mathbf{A} = (\alpha_{ir\ell\kappa})_{\substack{i=1,\ldots,m \\ r=1,\ldots,n \\ \ell=1,\ldots,2n \\ \kappa=1,\ldots,p}} \qquad \mathbf{X} = (x_{r\eta})_{\substack{r=1,\ldots,n \\ \eta=1,\ldots,4}}$$

$$h = (h_1, \ldots, h_p)^T \qquad G = (g_{ij})_{\substack{i=1,\ldots,m \\ j=1,\ldots,k}}.$$

Next, we relate $\Phi_\Sigma$ and $H$ and we write

$$Ce^{A(1-t)} B = \begin{bmatrix} \gamma_{11} & \cdots & \gamma_{1m} \\ \vdots & & \vdots \\ \gamma_{p1} & \cdots & \gamma_{pm} \end{bmatrix}.$$

We list the eigenvalues of $A$ as $a_r + ib_r$ for $r = 1, \ldots, n$ and let $\xi_{r\ell}(t) = \xi_\ell(a_r, b_r, t)$ for $r = 1, \ldots, n$ and $\ell = 1, \ldots, 2n$. Then there exists some $(\alpha_{ir\ell\kappa})$ such that

$$\gamma_{\kappa i}(t) = \sum_{r=1}^n \sum_{\ell=1}^{2n} \alpha_{ir\ell\kappa} \xi_{r\kappa}(t). \tag{5.4}$$

Let $\lambda = (\mathbf{A}, \mathbf{X}, h)$, where $\mathbf{A}$ satisfies (5.4), $\mathbf{X} = (x_{r\eta})$, where $x_{r1} = a_r, x_{r2} = b_r, x_{r3} = e^{a_r} \cos b_r$ and $x_{r4} = e^{a_r} \sin b_r$, and $h = Ce^A x^0$. We claim that

$$H(\lambda, G) = y(1) = \Phi_\Sigma(G\omega) \text{ for all } G \in \mathbb{R}^{mk}.$$

Note that the $\kappa$-th component of $\Phi_\Sigma(G\omega)$ is given by

$$\int_0^1 \sum_{i=1}^m \gamma_{\kappa i}(t) u_i(t) dt + h_\kappa$$

$$= \sum_{i=1}^m \sum_{r=1}^n \sum_{\ell=1}^{2n} \alpha_{ir\ell\kappa} \int_0^1 \xi_{r\ell}(t) \sum_{j=1}^k g_{ij} \omega_j(t) dt + h_\kappa$$

$$= \sum_{i=1}^m \sum_{r=1}^n \sum_{\ell=1}^{2n} \alpha_{ir\ell\kappa} \sum_{j=1}^k g_{ij} \int_0^1 \xi_\ell(a_r, b_r, t) \omega_j(t) dt + h_\kappa$$

$$= \sum_{i=1}^m \sum_{r=1}^n \sum_{\ell=1}^{2n} \alpha_{ir\ell\kappa} \sum_{j=1}^k g_{ij} \frac{\hat{P}_{\ell j}(x_{r1}, x_{r2}, x_{r3}, x_{r4})}{\hat{Q}_{\ell j}(x_{r1}, x_{r2}, x_{r3}, x_{r4})} + h_\kappa$$

$$= H_\kappa(\mathbf{A}, \mathbf{X}, h, G).$$

$\square$

Next we take $p = 1$ and study the VC-dimension of the sign system concept class, $\mathcal{C}_{m,1}$, where each control parameterized by $G$ gives rise to $\text{sign } y(1)$.

**Lemma 5.2.** *Sign system concept class $\mathcal{C}_{m,1}$ with initial condition $x(0) = 0$ satisfies*

$$VC(\mathcal{C}_{m,1}) \leq 2\left(2mn^2 + 4n\right) \log_2\left[8e\bigl(8mn^2 k(n + \ell_{\max}) + 1\bigr)\bigl(2nk + 2(1 + 2k)^n\bigr)\right],$$

*where $\ell_{\max}$ is given by (2.3).*

*Proof.* By Proposition 5.1 $y(1) = H(\lambda, G)$, where $\lambda \in \mathbb{R}^{2mn^2} \times \mathbb{R}^{4n}$ are considered as parameters. In fact, $y(1) = \frac{\mathbf{P}}{\mathbf{Q}}$, where $\mathbf{P}$ and $\mathbf{Q}$ denote piecewise polynomial functions. As in the statement of Goldberg-Jerrum bounds we have a function $F : \Lambda \times \mathbb{R}^{mk} \to \{0, 1\}$ defined by $F(\lambda, G) = \text{sign } H(\lambda, G)$. The concept class associated to the system identification problem is $\mathcal{F} := \{F(\lambda, \cdot) : \mathbb{R}^{mk} \to \{0, 1\} \,;\, \lambda \in \Lambda\}$, where $\Lambda = \mathbb{R}^{2mn^2 + 4n}$. Before applying the Goldberg-Jerrum bound we need to determine the possible degrees of $\mathbf{P}$ and $\mathbf{Q}$ with respect to the parameters.

The rationality condition implies that

$$\max_{\substack{i \leq \ell \leq n \\ 1 \leq j \leq k}} \left\{\deg(\hat{P}_{\ell j}), \deg(\hat{Q}_{\ell j})\right\} \leq d_{\max}.$$

Then

$$\frac{\widetilde{P}_{i\ell}}{\widetilde{Q}_{i\ell}} = \sum_{j=1}^k g_{ij} \frac{\hat{P}_{\ell j}}{\hat{Q}_{\ell j}}$$

so $\deg(\widetilde{Q}_{i\ell}) \leq k d_{\max}$ and $\deg(\widetilde{P}_{i\ell}) \leq k d_{\max}$. Note here that we are calculating the degree with respect to the system parameters, and the inputs $g_{ij}$ do not contribute.

By continuing in a similar fashion and combining $r$-summation to the $\ell$-summation in Proposition 5.1, we write $P_i/Q_i = \sum_{\ell=1}^{2n^2} \alpha_{i\ell} \widetilde{P}_{i\ell}/\widetilde{Q}_{i\ell}$ to conclude that $\deg(Q_i) \leq 2n^2 \deg(\widetilde{Q}_{i\ell}) = 2n^2 k d_{\max}$ and $\deg(P_i) \leq 2n^2 k d_{\max} + 1$. Finally, $\mathbf{P}/\mathbf{Q} = \sum_{i=1}^{m} P_i/Q_i$ with $\deg(\mathbf{Q}) \leq m2n^2 k d_{\max}$ and $\deg(\mathbf{P}) \leq m2n^2 k d_{\max} + 1$.

Recall that with $p = 1$ and initial condition $x(0) = 0$, using the notation of Proposition 5.1

$$y(1) = \sum_{i=1}^{m} \sum_{r=1}^{n} \sum_{\ell=1}^{2n} \alpha_{ir\ell} \sum_{j=1}^{k} g_{ij} \int_{0}^{1} \xi_\ell(x_{r1}, x_{r2}, t) \omega_j(t) dt.$$

The proof of Lemma 4.3 indicates that the denominator of

$$\int_{0}^{1} \xi_\ell(x_{r1}, x_{r2}, t) \omega_j(t) dt$$

equals

$$\left((x_{r1} + \alpha_j)^2 + (x_{r2} + \beta_j)^2\right)^{z_{\ell j}} \left((x_{r1} + \alpha_j)^2 + (x_{r2} - \beta_j)^2\right)^{z_{\ell j}},$$

where $\alpha_j, \beta_j$ are fixed parameters of the basis input function $\omega_j$ and $z_{\ell j} \in \mathbb{N}$.

By carrying out the summations we get $y(1) = \mathbf{P}/\mathbf{Q}$, where $\mathbf{Q}$ consists of powers of polynomials $f_{ij1}, f_{ij2}$ with

$$f_{ij1}(\mathbf{A}, \mathbf{X}, G) = (x_{i1} + \alpha_j)^2 + (x_{i2} + \beta_j)^2,$$
$$f_{ij2}(\mathbf{A}, \mathbf{X}, G) = (x_{i1} + \alpha_j)^2 + (x_{i2} - \beta_j)^2,$$

and $i = 1, \ldots, n$, $j = 1, \ldots, k$.

Our final step before applying the Goldberg-Jerrum bound is finding out the number of polynomial inequalities $s$ needed in the Boolean formula evaluating the sign of the final state output. This is done by studying the number of different $\mathbf{P}/\mathbf{Q}$ expressions without poles.

An upper bound for different $\mathbf{P}/\mathbf{Q}$ expressions without poles can be obtained by applying Theorem 4.6 to $2nk$ polynomials $f_{ij1}, f_{ij2}, i = 1, \ldots, n$ and $j = 1, \ldots, k$ viewing those as polynomials of $2n$ variables and each polynomial having degree 2. This gives the upper bound $(16ek)^{2n}$.

However, a more specific bound can be obtained in this problem. Note that varying $x_{i1}$ and $x_{i2}$ we can make at most one of the $2k$ polynomials $f_{ij1}, f_{ij2}, j = 1, \ldots, k$ to be zero. For example, $\gamma$ zeros among $f_{ij1}, f_{ij2}, i = 1, \ldots, n$ and $j = 1, \ldots, k$ can be obtained in $(2k)^\gamma \binom{n}{\gamma}$ ways and the number of possible sign assignments is obtained by summing over $\gamma$ yielding

$$\sum_{\gamma=0}^{n} (2k)^\gamma \binom{n}{\gamma} = (1 + 2k)^n.$$

Thus the number of $\mathbf{P}/\mathbf{Q}$ expressions without poles is $(1 + 2k)^n$, which gives rise to $2(1 + 2k)^n$ polynomials.

Note that in order to write $\operatorname{sign} y(1)$ as a Boolean formula evaluating polynomial inequalities and equalities it also has to include the $2nk$ polynomials $f_{ij1}, f_{ij2}, i =$ $1, \ldots, n$, $j = 1, \ldots, k$. Values of these polynomials determine which $\mathbf{P}/\mathbf{Q}$ expression is the valid one to determine sign $y(1)$. The Boolean formula for sign $y(1)$ can be given as a truth table involving polynomial inequalities of $2nk$ $f_{ij1}$, $f_{ij2}$ expressions and $2(1+2k)^n$ different $\mathbf{P}$ and $\mathbf{Q}$ expressions.

Using Lemma 4.3 for bound on $d_{\max}$, we apply the Goldberg-Jerrum bound with $s = 2nk + 2(2k+1)^n$, $d = m2n^2k4(n+\ell_{\max}) + 1$ and $\ell = 2mn^2 + 4n$. □

A simple example of a piecewise polynomial function $\mathbf{P}/\mathbf{Q}$ together with the decision table for the final output is provided in Appendix A.

*Remark 5.3.* The VC-dimension bound is modified for the more abstract rationality conditions as follows. Evaluating the sign of the output involves the evaluation of $2(8ed_{\max}2n^2kh/4n)^{4n} + 2n^2kh$ polynomials; $2n^2kh$ evaluations are needed to find an appropriate piece and by Theorem 4.6 the maximum number of possible expressions of the type $\mathbf{P}/\mathbf{Q}$ is bounded by $(8ed_{\max}2n^2kh/4n)^{4n}$. Applying the Goldberg-Jerrum bound with $s = 2(8ed_{\max}2n^2kh/4n)^{4n} + 2n^2kh$, $d = m2n^2kd_{\max} + 1$ and $\ell = 2mn^2 + 4n$ gives the result.

**Proof of Theorem 3.1, the VC-dimension upper bound, p = 1.** By using the previous notation $y = Ce^A x^0 + C \int_0^1 e^{A(1-t)} Bu(t) dt$. Let $\widetilde{x} = Ce^A x^0$. Then $y = \widetilde{x} + \mathbf{P}/\mathbf{Q} = (\widetilde{x}\mathbf{Q} + \mathbf{P})/\mathbf{Q} = \widetilde{\mathbf{P}}/\mathbf{Q}$. This has $2mn^2 + 4n + 1$ parameters and $\deg(\widetilde{\mathbf{P}}) \leq m2n^2kd_{\max} + 1$. □

*Remark 5.4.* Taking $2n^2$ functions $\xi_{r\ell}(t)$, $r = 1, \ldots, n$, $\ell = 1, \ldots, 2n$ in Proposition 5.1 is clear overcounting. Further calculations indicate that $O(n \ln n)$ functions $\xi_{r\ell}(t)$ would be enough. For more details see Appendix B.

## 5.2 Lower Bounds for the VC-Dimension

The lower bounds for the VC-dimension are developed for a single-input single-output system with initial state zero. The control is

$$u = \sum_{j=1}^{k} g_j \omega_j : [0,1] \to \mathbb{R}.$$

We derive lower bounds by fixing the structure of $A$, $B$, and $C$, and using the dual VC-dimension and axis shattering following the ideas of DasGupta and Sontag [7]. Lemmata 5.6, 5.7 and 5.8 given in this section together prove Theorem 3.2. These lower bounds are very general; we just assume that the input functions are continuous and linearly independent, thus no particular structure on input function is required as in the upper bounds.

To make the next proof cleaner we formulate a part of it as a separate proposition. (The proposition is a standard fact but we include a short proof for the completeness of exposition.)

**Proposition 5.5.** *Let $\omega_j : [0,1] \to \mathbb{R}$, $j = 1, \ldots, k$ be continuous and linearly independent. Then the functions*

$$h_j(\lambda) = \int_0^1 e^{\lambda t} \omega_j(t) dt, \quad j = 1, \ldots, k$$

*are linearly independent.*

*Proof.* We begin by showing that for any $f \in C[0,1] = \{f \colon [0,1] \to \mathbb{R} \,;\, f \text{ continuous}\}$

$$\int_0^1 e^{\lambda t} f(t) dt = 0 \quad \text{for all } \lambda \in \mathbb{R}$$

implies $f(t) \equiv 0$ on $[0,1]$. The above integral is an analytic function of $\lambda$. Moreover, all derivatives in $\lambda$ are zero at $\lambda = 0$, which gives that $\int_0^1 t^k f(t) dt = 0$ for all $k = 1, 2, \ldots$. As $f$ is continuous and polynomials are dense in $C[0,1]$, we have that $f(t) \equiv 0$.

Linear independence of $h_1, \ldots, h_k$ follows: If $\alpha_1 h_1(\lambda) + \cdots + \alpha_k h_k(\lambda) = 0$ for all $\lambda \in \mathbb{R}$, then

$$\int_0^1 e^{\lambda t}(\alpha_1 \omega_1(t) + \cdots + \alpha_k \omega_k(t)) dt = 0, \quad \forall \lambda \in \mathbb{R}$$

and by the above argument $\alpha_1 \omega_1(t) + \cdots + \alpha_k \omega_k(t) = 0$ for all $t$. As $\omega_1, \ldots, \omega_k$ were assumed to be linearly independent, we have that $(\alpha_1, \ldots, \alpha_k) = (0, \ldots, 0)$. $\square$

**Lemma 5.6 (Lower bound 1).** *Sign system concept class $\mathcal{C}_{1,1}$ with scalar inputs and scalar outputs satisfies*

$$VC(\mathcal{C}_{1,1}) \geq m' \left\lfloor \log_2 \left\lfloor \frac{k}{m'} \right\rfloor \right\rfloor,$$

*where $m' = \min\{n, k\}$.*

*Proof.* Let $\omega_j(t)$, $j = 1, \ldots, k$ be continuous and linearly independent. Let $A$ have $n$ distinct real eigenvalues $-\lambda_1, \ldots, -\lambda_n$, and take $B$ and $C$ so that

$$Ce^{A(1-t)}B = \sum_{i=1}^{m'} e^{\lambda_i t},$$

where $m' = \min\{n, k\}$. Then the final output of the system is

$$y(1) = \int_0^1 Ce^{A(1-t)}B \sum_{j=1}^k g_j \omega_j(t) dt = \sum_{i=1}^{m'} \sum_{j=1}^k g_j \int_0^1 e^{\lambda_i t} \omega_j(t) dt.$$

Define $h_j(\lambda) = \int_0^1 e^{\lambda t} \omega_j(t) dt$. By Proposition 5.5 the $h_j$'s are linearly independent and we can find $\lambda_1, \ldots, \lambda_k$ such that the matrix

$$\begin{bmatrix} h_1(\lambda_1) & \cdots & h_k(\lambda_1) \\ \vdots & & \vdots \\ h_1(\lambda_k) & \cdots & h_k(\lambda_k) \end{bmatrix}$$

has rank $k$.

The control system with sign observations gives the mapping $F : \mathbb{R}^{m'} \times \mathbb{R}^k \to \{0, 1\}$ by

$$(\lambda_1, \ldots, \lambda_{m'}, g_1, \ldots, g_k) \mapsto \operatorname{sign}\left[\sum_{i=1}^{m'} \sum_{j=1}^{k} g_j h_j(\lambda_i)\right].$$

We show that the mapping from parameters $\lambda_1, \ldots, \lambda_{m'}$ to $\{0, 1\}$ can be axis shattered. Let $L = \{\lambda_1, \ldots, \lambda_k\}$ be so that $[h_j(\lambda_i)]_{i,j}$ has rank $k$. Denote by $L_1, \ldots, L_{m'}$ disjoint subsets of $L$ such that $|L_i| = \lfloor k/m' \rfloor$ and let $M = L \setminus \{\bigcup_{i=1}^{m'} L_i\}$. Next we want to interpolate in the points of $L$.

Fix $s$, $1 \leq s \leq m'$ and let $\phi : L_s \to \{0, 1\}$ be any dichotomy. Next find $g_1, \ldots, g_k$ such that

$$\begin{aligned}
\sum_{j=1}^{k} g_j h_j(\lambda_s) &= \phi(\lambda_s), && \forall \lambda_s \in L_s, \\
\sum_{j=1}^{k} g_j h_j(\lambda) &= 0, && \forall \lambda \in (L \cup M) \setminus L_s.
\end{aligned} \tag{5.5}$$

Let $g_1^*, \ldots g_k^*$ satisfy (5.5). (A unique solution exists because $[h_j(\lambda_i)]$ has rank $k$.) Then

$$F[\lambda_1, \ldots, \lambda_{m'}, g_1^*, \ldots, g_k^*] = \operatorname{sign}\left[\sum_{i=1}^{m'} \sum_{j=1}^{k} g_j^* h_j(\lambda_i)\right] = \phi(\lambda),$$

when $\lambda \in L_s$ and for all $(\lambda_1, \ldots, \lambda_{m'}) \in L_1 \times \cdots \times L_{m'}$.

Let $\widetilde{\mathcal{F}} = \{F(\lambda_1, \ldots, \lambda_{m'}, \cdot) : \mathbb{R}^k \to \{0, 1\} \,;\, (\lambda_1, \ldots, \lambda_{m'}) \in \mathbb{R}^{m'}\}$. By the Axis shattering bound given in Theorem 4.4

$$\operatorname{VC}(\widetilde{\mathcal{F}}) \geq m' \left\lfloor \log_2 \left\lfloor \frac{k}{m'} \right\rfloor \right\rfloor$$

and thus $\operatorname{VC}(\mathcal{C}_{1,1}) \geq \operatorname{VC}(\widetilde{\mathcal{F}})$, where $\mathcal{C}_{1,1}$ is the control system concept class with $p = m = 1$. $\square$

**Lemma 5.7 (Lower bound 2).** *If $k \leq n$ then*

$$VC(\mathcal{C}_{1,1}) \geq k.$$

*Proof.* We make a small modification of the above argument. Assume that $k \leq n$ and let $A$ have $n$ real eigenvalues $\lambda_1, \ldots, \lambda_n$. Next we take $B$ and $C$ so that $Ce^{A(1-t)}B = \sum_{i=1}^{n} e^{\lambda_i t} \beta_i$, where $(\beta_1, \ldots, \beta_n, \lambda_1, \ldots, \lambda_n)$ are considered as system parameters.

We study the mapping

$$(\beta_1, \ldots, \beta_n, \lambda_1, \ldots, \lambda_n, g_1, \ldots, g_k) \mapsto \operatorname{sign}\left[\sum_{i=1}^{n} \sum_{j=1}^{k} g_j h_j(\lambda_i) \beta_i\right]$$

$$= \operatorname{sign}\left[\sum_{j=1}^{k} g_j \underbrace{\sum_{i=1}^{n} h_j(\lambda_i) \beta_i}_{\gamma_j}\right] = \operatorname{sign}\left[\sum_{j=1}^{k} g_j \gamma_j\right].$$

Given $(\gamma_1, \ldots, \gamma_k)$, by linear independence of $h_1, \ldots, h_k$, we can find $\lambda_1, \ldots, \lambda_n, \beta_1, \ldots, \beta_n$ such that $\sum_{i=1}^n h_j(\lambda_i)\beta_i = \gamma_j$, $j = 1, \ldots, k$. But $(\gamma_1, \ldots, \gamma_k)$ can be viewed as a normal vector for a hyperplane through the origin in $\mathbb{R}^k$ and the concept class associated to the mapping $(g_1, \ldots, g_k) \mapsto \text{sign}\left[\sum_{j=1}^k g_j \gamma_j\right]$ as $(\gamma_1, \ldots, \gamma_k)$ varies has VC-dimension $k$. Hence $\text{VC}(\mathcal{C}_{1,1}) \geq k$. □

**Lemma 5.8 (Lower bound 3).** *If $n \leq k$ then*
$$VC(\mathcal{C}_{1,1}) \geq n.$$

*Proof.* Our construction for the control system is as in the previous proof but now we assume that $n \leq k$ and we study

$$(\beta_1, \ldots, \beta_n, \lambda_1, \ldots, \lambda_n, g_1, \ldots, g_k) \mapsto \text{sign}\left[\sum_{i=1}^n \sum_{j=1}^k g_j h_j(\lambda_i)\beta_i\right]$$

$$= \text{sign}\left[\sum_{i=1}^n \underbrace{\sum_{j=1}^k g_j h_j(\lambda_i)}_{\tilde{g}_i} \beta_i\right] = \text{sign}\left[\sum_{i=1}^n \tilde{g}_i \beta_i\right]$$

and again by linear independence and the above hyperplane argument (now via first transforming $(g_1, \ldots, g_k)$) we can conclude that the above mapping has VC-dimension $n$. Thus $\text{VC}(\mathcal{C}_{1,1}) \geq n$. □

## 5.3 VC-Dimension Upper Bounds for *p*-Dimensional Outputs

We begin by proving Theorem 3.3.

**Proof of the VC-dimension upper bound.** We develop an upper bound based on the bound for a scalar sign-observation. We have seen that under the rationality assumption (RAT) the scalar output is a piecewise rational expression $\mathbf{P}/\mathbf{Q}$. In general, the control system maps $G$ to $(\text{sign}(\mathbf{P}_1/\mathbf{Q}_1), \ldots, \text{sign}(\mathbf{P}_p/\mathbf{Q}_p))^T$, which is understood as a binary representation of a number in $\{0, 1, \ldots, 2^p - 1\}$. Let $f : \mathbb{R}^{mk} \to \{0, \ldots, 2^p - 1\}$ be the mapping given by the control system, and denote the class of all such mappings by $\mathcal{F}$. For each $f \in \mathcal{F}$ introduce a loss function $L_{0\text{-}1,f}(z, a) = L_{0\text{-}1}(f(z), a) = 1$, when $f(z) \neq a$, and 0 otherwise. Define the class $L_{0\text{-}1,\mathcal{F}} = \{L_{0\text{-}1,f} \ ; \ f \in \mathcal{F}\}$.

In order to calculate the value of the output, after determining an appropriate piece, one needs to know the truth values of the expressions $\mathbf{P}_1 > 0, \mathbf{Q}_1 > 0, \ldots, \mathbf{P}_p > 0$ and $\mathbf{Q}_p > 0$, where $\mathbf{P}$'s and $\mathbf{Q}$'s are polynomials on inputs and parameters of the control system. To evaluate the value of the loss function $L_{0\text{-}1,f}(z, a)$, one needs the truth values of $y = 0, y = 1, \ldots, y = 2^p - 2$.

In the general case one needs $2nk + 2p(2k+1)^n + 2^p - 1$ truth values. As this procedure evaluates only polynomials, we can use the Goldberg-Jerrum bound again. The maximum degree of the polynomials is $m2n^2k4(n + \ell_{\max}) + 1$, and the total number of parameters is $2pn^2m + 4n + p$, where the last term comes from the initial condition. □

*Remark 5.9.* An upper bound can be derived by using the uniform $\Psi$-dimension defined by Ben-David et. al. [6]. As before let $f : \mathbb{R}^{mk} \to \{0, \ldots, 2^p - 1\}$ be the mapping given by the control system. Let $\Psi = \{\psi^1, \ldots, \psi^p\}$ be a collection of distinguishers $\psi^j : \{0, \ldots, 2^p - 1\} \to \{0, 1\}$ given by $\psi^j(\text{sign}(y_1), \ldots, \text{sign}(y_p)) = \text{sign}(y_j)$.

Let $\psi_f^j : \mathbb{R}^{mk} \to \{0, 1\}$ be given by $\psi_f^j(z) = \psi^j(f(z))$ and $\psi_\mathcal{F}^j = \{\psi_f^j \ ; \ f \in \mathcal{F}\}$. An upper bound for $\text{VC}(\psi_\mathcal{F}^j)$ is given by Theorem 3.1. Ben-David et. al. define the uniform $\Psi$-dimension of $\mathcal{F}$, in which only one distinguisher is selected to shatter, as

$$\Psi - \dim_U(\mathcal{F}) = \max\{\text{VC}(\psi_\mathcal{F}^j) \ ; \ \psi^j \in \Psi\}.$$

That is, the upper bound for the uniform $\Psi$-dimension is just the upper bound obtained earlier for scalar observations. However, the uniform $\Psi$-dimension can not be compared to the VC-dimension as such. For example, the formula for sample complexity with $\Psi$-dimension is different: to calculate the sample complexity the uniform $\Psi$-dimension is multiplied by $p \log_2(2^p - 1)$, see [6, Theorem 27].

## 6 A Fat-Shattering Bound

We begin this section by proving Theorems 3.5 and 3.6. As a corollary of Theorem 3.6 we prove the fat-shattering bound appearing in Theorem 2.2 bounding the sample complexity for proper agnostic learning.

**Proof of Theorem 3.5.** For the first part of the proof we use a generic set $B$ for the parameters. Assume that we can $\gamma$-shatter a set of inputs $\{u_1, \ldots, u_d\}$ and there exists $\{r_1, \ldots, r_d\}$ such that, for each assignment $b \in \{0, 1\}^d$, there exists a $\lambda \in B$ such that

$$F(\lambda, u_i) \geq r_i + \gamma, \quad \text{if } b_i = 1, \text{ and}$$
$$F(\lambda, u_i) \leq r_i - \gamma, \quad \text{otherwise.}$$

We write $\lambda \sim \mu$ if and only if the parameters $\lambda$ and $\mu$ give the same assignment for all $\{u_1, \ldots, u_d\}$. Further, let $\Lambda = \{\lambda_1, \ldots, \lambda_{2^d}\}$ be a collection of parameters that shatter $\{u_1, \ldots, u_d\}$ and let $\lambda_i, \lambda_j \in \Lambda$. Now $\lambda_i \not\sim \lambda_j$ implies that there exists $u^* \in \{u_1, \ldots, u_d\}$ and $r^* \in \{r_1, \ldots, r_d\}$ such that $F(\lambda_i, u^*) \geq \gamma + r^*$ and $F(\lambda_j, u^*) \leq \gamma - r^*$, or vice versa. Hence

$$2\gamma \leq |F(\lambda_i, u^*) - F(\lambda_j, u^*)| \leq L\|\lambda_i - \lambda_j\|$$

and so $\|\lambda_i - \lambda_j\| \geq 2\gamma/L$. That is, the set $\Lambda$ of cardinality $2^d$ is a $2\gamma/L$-separated set in $B$. Now the fat-shattering bounds follow by calculating $2\gamma/L$-packing numbers for different sets $B$.

If $B = B_\infty^k(C)$ the maximum possible cardinality for an $\epsilon$-separated set is $\lfloor 2C/\epsilon \rfloor^k$, and thus

$$2^d \leq \left\lfloor \frac{2C}{2\gamma/L} \right\rfloor^k = \left\lfloor \frac{CL}{\gamma} \right\rfloor^k$$

and solving for $d$ yields $d \leq k \log_2 \lfloor CL/\gamma \rfloor$.

Similarly, if $B = \bar{B}_\infty^k(C)$ the maximum possible cardinality for an $\epsilon$-separated set is $(1 + \lfloor 2C/\epsilon \rfloor)^k$ and by a similar argument we arrive at the bound $d \leq k \log_2 (1 + \lfloor LC/\gamma \rfloor)$.

For $B = \bar{B}_2^k(C)$, let $P(\epsilon)$ be a collection of $\epsilon$-separated sets in $\bar{B}_2^k(C)$ and let $|P(\epsilon)|$ denote its cardinality. As all open balls with radius $\epsilon/2$ with centers at $\epsilon$-separated points have to be disjoint and their union has to be inside a ball of radius $C + \epsilon/2$, we get that $|P(\epsilon)|\alpha(k)(\epsilon/2)^k \leq \alpha(k)(C+\epsilon/2)^k$, where $\alpha(k) = \pi^{k/2}/\Gamma(k/2+1)$ is the volume of a unit ball in $\mathbb{R}^k$. Hence $|P(\epsilon)| \leq (C + 2/\epsilon)^k$ and

$$2^d \leq (C + L/\gamma)^k,$$

i.e., $d \leq k \log_2(C + L/\gamma)$. $\square$

Next we prove Theorem 3.6 by applying the Lipschitz bound to a control system.

**Proof of Theorem 3.6.** Our aim is to compute the Lipschitz constant associated to the control system in Definition 2 and then we apply Theorem 3.5.

Denote the system parameters $(\alpha_{11}, \ldots, \alpha_{nm}, a_1, b_1, \ldots, a_r, b_r)$ by $\lambda$ and assume $\|\lambda\|_\infty < 1$. Let

$$F(\lambda, u) = y(\tau) = \int_0^\tau \sum_{i=1}^m \sum_{\ell=1}^n \alpha_{i\ell} \xi_\ell(t) u_i(\tau - t) dt.$$

Functions $\xi_1(t), \ldots, \xi_n(t)$ are of the form $\xi(t) = t^c e^{at} \sin(bt)$ or $\xi(t) = t^c e^{at} \cos(bt)$ where $a + ib$ is an eigenvalue of $A$ and $c \in \{0, \ldots, n-1\}$. Thus taking a partial derivative with respect to $a$ or $b$ will increase the power of $t$ by one and change the trigonometric functions. Therefore,

$$\left|\frac{\partial F(\lambda, u)}{\partial \alpha_{\kappa\rho}}\right| = \left|\frac{\partial y(\tau)}{\partial \alpha_{\kappa\rho}}\right| = \left|\sum_{i=1}^m \sum_{\ell=1}^n d(i, \ell) \int_0^\tau \xi_\ell(t) u_i(\tau - t) dt\right|$$

$$\leq nm \int_0^\tau |\xi_\ell(t) u_i(\tau - t)| \, dt \leq nm\tau^{n-1} e^\tau M,$$

because $\sup_{t \in [0,\tau]} |\xi_\ell(t)| \leq e^\tau \tau^n$ and $d(i, \ell) = \partial \alpha_{ij}/\partial \alpha_{\kappa\rho} = 1$ if $(i, \ell) = (\kappa, \rho)$ and zero otherwise. Similarly we calculate

$$\left|\frac{\partial F(\lambda, u)}{\partial a_\kappa}\right| = \left|\frac{\partial y(\tau)}{\partial a_\kappa}\right| \leq nm\tau^n e^\tau M \text{ and}$$

$$\left|\frac{\partial F(\lambda, u)}{\partial b_\kappa}\right| = \left|\frac{\partial y(\tau)}{\partial b_\kappa}\right| \leq nm\tau^n e^\tau M$$

as $\sup_{t \in [0,\tau]} \left|\frac{\partial \xi_\ell(t)}{\partial_\kappa}\right| \leq e^\tau \tau^n$ and $\sup_{t \in [0,\tau]} \left|\frac{\partial \xi_\ell(t)}{\partial b_\kappa}\right| \leq e^\tau \tau^n$.

Now the Lipschitz constant can be taken to be $L = n^2 m e^\tau \tau^n M$ as

$$|F(\lambda, u) - F(\lambda^*, u)| = |\nabla F \cdot (\lambda - \lambda^*)| \leq L\|\lambda - \lambda^*\|_\infty.$$

The number of system parameters is at most $nm + n = (m+1)n$ and we get the level fat-shattering bound by applying Theorem 3.5 with space dimension $n(m+1)$ and $L = n^2 m e^\tau \tau^n M$. $\square$

As a corollary, we combine the above result together with a pseudo-dimension bound to prove the fat-shattering bound given in Theorem 2.2.

**Corollary 6.1 (Fat-shattering bound in Theorem 2.2).** *Assume that the system $\dot{x} = Ax + Bu$, $y = Cx$, $x(0) = 0$ can be parameterized by $\lambda \in \mathbb{R}^{n(m+1)}$ as in Definition 2 with $\|\lambda\|_\infty < 1$ and assume in addition that the control is given by $u = G\omega$ where the input basis functions $\omega_j$ are in $\Omega$ given by (2.2). We denote the corresponding control system class by*

$$\mathcal{F}_B = \{F(\lambda, \cdot) : U \to \mathbb{R} \;;\; \lambda \in B\}.$$

*Then*

$$\mathrm{fat}_\gamma(\mathcal{F}_B) \leq \min \begin{cases} (m+1)n \log_2 \left\lfloor \frac{n^2 m \tau^n e^\tau kM}{\gamma} \right\rfloor, \\ 2(m+4)n \log_2 \left( 8e\bigl(nmk4(n+\ell_{\max})+1\bigr)\bigl(2nk + 2(2k+1)^n\bigr) \right), \end{cases}$$

*where $\ell_{\max}$ is given by (2.2) and (2.3) and $M$ is a constant satisfying*

$$\int_0^\tau |u_i(\tau - t)| dt \leq kM$$

*for all $i = 1, \ldots, m$.*

*Proof.* The first part of the bound follows from Theorem 3.6 with $kM$ in place of $M$.

The remaining part of the bound comes from the pseudo-dimension bound. First we derive the associated VC-dimension bound. As we assumed that $A$ has a fixed Jordan block structure, every entry of $e^{A(1-t)}$ is a linear combination of $n$ functions $\xi_1(t), \ldots, \xi_n(t)$. (That is, we don't need to consider all possible functions over different Jordan block structures.) This implies that in the Goldberg-Jerrum argument of Section 5.1 we can take $\ell = mn + 4n$, $d = nmk4(n+\ell_{\max}) + 1$ and $s = 2nk + 2(2k+1)^n$. Moreover, in that section the VC-dimension bounds were derived for time interval $[0, 1]$. However, the upper bound depends on the number of system parameters and the degrees of polynomials to be evaluated. Changing the time interval to be $[0, \tau]$ means just that we replace the eigenvalue parameters (referring to the proof of Proposition 5.1) $a, b, e^a \cos b, e^a \sin b$ by $a\tau, b\tau, e^{a\tau} \cos b\tau, e^{a\tau} \sin b\tau$.

The above bound is also a bound for the pseudo-dimension. Observe that for $\mathcal{G} = \{g : X \to \mathbb{R}\}$, the pseudo-dimension can be defined as $\mathrm{PD}(\mathcal{G}) = \mathrm{VC}\{\mathrm{Ind}(x, y) = \mathrm{sign}(g(x) - y) \;;\; g \in \mathcal{G}\}$. Hence we study the VC-dimension associated to $\mathrm{sign}(y(\tau) - z) = \mathrm{sign}(\mathbf{P}/\mathbf{Q} - z) = \mathrm{sign}(\hat{\mathbf{P}}/\mathbf{Q})$, where $\hat{\mathbf{P}} = \mathbf{P} - z\mathbf{Q}$ has the same degree as $\mathbf{P}$ with respect to the parameters. Here $z$ is a new input, but the bound utilizing Goldberg-Jerrum technique does not depend on the dimension of the inputs, and hence the above VC-dimension bound is also a bound for the pseudo-dimension. (Note that here in the scale sensitive setting we do not apply the pseudo-dimension results of Section 5.3 using loss functions, as those rescaled the outputs.) □

*Remark 6.2.* In the special case of scalar controls (i.e., $m = 1$) that are given by a linear combination of $k$ fixed input functions, the control system class discussed can be viewed as a family of linear mappings in $\mathbb{R}^k$ as calculated in (2.5). Let $\mathcal{F} = \{x \in \mathbb{R}^k \mapsto \sum_i w_i x_i + \theta \; ; \; \|w\|_2 = 1\}$. If we further restrict $x$ such that $\|x\|_2 \leq R$ and $|\theta| \leq R$ then [18, 16] have shown that

$$\text{fat}_\gamma(\mathcal{F}) \leq \min\{9R^2/\gamma^2, k+1\} + 1.$$

In this special case our $\gamma$-shattering result is of the form $c_1 \log_2(\lfloor c_2 k/\gamma \rfloor)$, where $c_1$ and $c_2$ are constants. This gives improvement over the hyperplane bound when the margin $\gamma$ is small.

# 7 A Class of Systems with VC-Dimension $k$

For the control system (2.1) with scalar control $u(t) = \sum_{i=1}^{k} g_i \omega_i(t)$ and unrestricted $\omega_1, \ldots, \omega_k$, the standard half-space argument gives an upper bound $k$. This bound is tight. We will give an example of a single-input, single-output one parameter family of control systems in dimension two that has VC-dimension $k$, when the controls are of the form $u(t) = \sum_{i=1}^{k} g_i \omega_i(t)$ and $\omega_i(1-t) = 1_{[2^{-i}, 2^{-i}+2^\alpha]}$, where $\alpha = -2(k+1)$.

Consider a control system

$$\begin{aligned} \dot{x}_1 &= x_2 \\ \dot{x}_2 &= -\lambda^2 x_1 + u \\ y &= -x_1. \end{aligned} \qquad (7.1)$$

For time interval [0,1] and initial condition $(x_1, x_2) = (0,0)$, the output is given by $y(1) = \int_0^1 \sin(\lambda t) u(1-t) dt$.

**Lemma 7.1.** *Controls $\{\omega_1, \ldots, \omega_k\}$ such that $\omega_i(1-t) = 1_{[2^{-i}, 2^{-i}+2^\alpha]}$, where $\alpha = -2(k+1)$, are shattered by the control system (7.1) with sign-observations.*

*Proof.* Let $T = \{2^{-i}, i = 1, \ldots, k\}$ and $J \subseteq T$. Define $\lambda_J = \pi \sum_{i=1}^{k} a_i 2^i$, where $a_i = 1$ if $2^{-i} \notin J$ and $a_i = 0$ otherwise. Now if $t = 2^{-\ell}$, then

$$\lambda_J t = \pi \sum_{i=1}^{k} a_i 2^{i-\ell} = \pi \left( \underbrace{\sum_{i=1}^{\ell-1} a_i 2^{i-\ell}}_{1/2\, c_2} + a_\ell + \underbrace{\sum_{i=\ell+1}^{k} a_i 2^{i-\ell}}_{2c_1} \right),$$

where $c_1 \in \mathbb{N}$ and $0 \leq c_2 < 2$.

Hence $\sin(\lambda_J t) = \sin(\pi(1/2\, c_2 + a_\ell))$. Note that if $a_\ell = 0$ then $1/2\, c_2 + a_\ell \in [0, 1 - 2^{-\ell}]$ and if $a_\ell = 1$ then $1/2\, c_2 + a_\ell \in [1, 2 - 2^{-\ell}]$. Thus $\sin(\pi(1/2\, c_2 + a_\ell)) \geq 0$ if $a_\ell = 0$ and $\sin(\pi(1/2\, c_2 + a_\ell)) \leq 0$ if $a_\ell = 1$. Therefore,

$$\sin(\lambda_J t) \geq 0 \Leftrightarrow a_\ell = 0$$

and

$$\sin(\lambda_J t) \geq 0 \Leftrightarrow t \in J.$$

Further,
$$\int_{2^{-\ell}}^{2^{-\ell}+2^\alpha} \sin(\lambda_J t)dt \geq 0 \Leftrightarrow a_\ell = 0,$$

where $\alpha$ is taken so that
$$\sum_{j=1}^{k} 2^j 2^\alpha \leq 2^{-(k+1)}. \tag{7.2}$$

This assures that when $\ell \leq k$ and $t \in [2^{-\ell}, 2^{-\ell} + 2^\alpha]$
$$\lambda_J t \in [0, \pi(1 - 2^{-\ell} + \sum_{j=1}^{k} 2^j 2^\alpha)] \subset [0, \pi), \text{ if } a_\ell = 0$$

or similarly
$$\lambda_J t \in [\pi, 2\pi), \text{ if } a_\ell = 1.$$

In (7.2) we can take $\alpha = -2(k+1)$ as $\sum_{j=1}^{k} 2^j = 2^{k+1} - 2$.

In this way the integrand in
$$\int_{2^{-\ell}}^{2^{-\ell}+2^\alpha} \sin(\lambda_J t)dt$$

is either positive or negative.

For $S \subseteq \{1, \ldots, k\}$, let $J = \{2^{-i}, i \in S\}$. For each $\omega_i$,
$$\int_0^1 \sin(\lambda_J t)\omega_i(1-t)dt = \int_{2^{-i}}^{2^{-i}+2^\alpha} \sin(\lambda_J t)dt > 0 \Leftrightarrow i \in S,$$

i.e., the set of controls $\{\omega_1, \ldots, \omega_k\}$ is shattered by the mapping
$$\omega_i \mapsto \text{sign}\left[\int_0^1 \sin(\lambda_J t)\omega_i(1-t)dt\right].$$

$\square$

# A  An Example of the Goldberg-Jerrum Bound

We begin this appendix with an informal discussion on the Goldberg-Jerrum technique used to prove the VC-dimension upper bounds in this paper.

We want to write $y(1) = \mathbf{P}/\mathbf{Q}$, where $\mathbf{P}$ and $\mathbf{Q}$ are polynomials. Unfortunately, the value of sign $y(1)$ can not be obtained by just evaluating $\mathbf{P}$ and $\mathbf{Q}$ since $\mathbf{Q}$ may have zeros. Therefore we need to write
$$y(1) = \begin{cases} \mathbf{P}_1/\mathbf{Q}_1, & \text{if } f_1 \neq 0, \ldots, f_\mu \neq 0 \\ \vdots \\ \mathbf{P}_\gamma/\mathbf{Q}_\gamma, & \text{if } f_1 = 0, \ldots, f_\mu \neq 0 \end{cases}$$

so that after evaluating $\mu$ polynomials $f_1, \ldots, f_\mu$ we can pick a definition $\mathbf{P}_i/\mathbf{Q}_i$ without poles in a region defined by the $\mu$ polynomials. When $y(1)$ is defined in this way $\operatorname{sign} y(1)$ can be easily expressed by a Boolean formula evaluating $2\gamma + \mu$ polynomial inequalities and equalities.

For simplicity we assume that $p = 1$ and the initial condition $x(0) = 0$. Then using the notation of Proposition 5.1 we write

$$y(1) = \sum_{i=1}^{m} \sum_{r=1}^{n} \sum_{\ell=1}^{2n} \alpha_{ir\ell} \sum_{j=1}^{k} g_{ij} \int_0^1 \xi_\ell(x_{r1}, x_{r2}, t)\omega_j(t)dt$$

and by the proof of Lemma 4.3

$$\int_0^1 \xi_\ell(x_{r1}, x_{r2}, t)\omega_j(t)dt = \frac{P_{\ell j}}{((x_{r1} + \alpha_j)^2 + (x_{r2} + \beta_j)^2)^{z_{\ell j}} ((x_{r1} + \alpha_j)^2 + (x_{r2} - \beta_j)^2)^{z_{\ell j}}},$$

where $P_{\ell j}$ is some polynomial, $z_{\ell j} \in \mathbb{N}$ and $\omega_j(t) = t^{\ell_j} e^{\alpha_j t} \sin(\beta_j t)$ or $\omega_j(t) = t^{\ell_j} e^{\alpha_j t} \cos(\beta_j t)$. Hence the denominator of

$$\sum_{j=1}^{k} g_{ij} \int_0^1 \xi_\ell(x_{r1}, x_{r2}, t)\omega_j(t)dt$$

is

$$\left((x_{r1} + \alpha_1)^2 + (x_{r2} + \beta_1)^2\right)^{z_{\ell 1}} \left((x_{r1} + \alpha_1)^2 + (x_{r2} - \beta_1)^2\right)^{z_{\ell 1}} \times$$
$$\cdots \times \left((x_{r1} + \alpha_k)^2 + (x_{r2} + \beta_k)^2\right)^{z_{\ell k}} \left((x_{r1} + \alpha_k)^2 + (x_{r2} - \beta_k)^2\right)^{z_{\ell k}}.$$

By carrying out all summations $y(1) = \mathbf{P}/\mathbf{Q}$. The denominator $\mathbf{Q}$ consists of the product

$$\prod_{r=1}^{n} \left(\left((x_{r1} + \alpha_1)^2 + (x_{r2} + \beta_1)^2\right)^* \left((x_{r1} + \alpha_1)^2 + (x_{r2} - \beta_1)^2\right)^* \times \right.$$
$$\left. \cdots \times \left((x_{r1} + \alpha_k)^2 + (x_{r2} + \beta_k)^2\right)^* \left((x_{r1} + \alpha_k)^2 + (x_{r2} - \beta_k)^2\right)^*\right),$$

where $*$'s stand for some unspecified powers. Hence the zeros of $\mathbf{Q}$ are determined by $2nk$ polynomials

$$f_{ij1} = (x_{i1} + \alpha_j)^2 + (x_{i2} + \beta_j)^2,$$
$$f_{ij2} = (x_{i1} + \alpha_j)^2 + (x_{i2} - \beta_j)^2,$$

and $i = 1, \ldots, n$, $j = 1, \ldots, k$. The number of different sign assignments determining $\gamma$ is calculated as in the proof of Lemma 5.2.

**Example**

The purpose of the following example is to illustrate the function $y = \mathbf{P}/\mathbf{Q}$ used in the Goldberg-Jerrum technique together with the sequence of polynomial evaluations involved and a table for the final output depending on the outcomes of the polynomial evaluations.

Take $m = 1$, $n = 2$, $k = 2$, and assume that $A$ has complex eigenvalues $a \pm ib$. Take basis input functions to be $\omega_1(t) = e^t$ and $\omega_2(t) = e^{2t}$. Then

$$y(1) = \sum_{l=1}^{2} \alpha_l \sum_{j=1}^{2} g_j \int_0^1 \xi_l(t)\omega_j(t)\,dt,$$

where $\xi_1(t) = e^{at}\sin(bt)$, $\xi_2(t) = e^{at}\cos(bt)$, $\alpha_1$, $\alpha_2$, $a$, $b$, $e^a \sin b$ and $e^a \cos b$ are system parameters and $g_1$, $g_2$ are input parameters.

By using formulas

$$\int_0^1 e^{\tilde{a}t} \sin(\tilde{b}t)\,dt = \frac{e^{\tilde{a}}(\tilde{a}\sin\tilde{b} - \tilde{b}\cos\tilde{b}) + \tilde{b}}{\tilde{a}^2 + \tilde{b}^2}, \quad \text{and}$$

$$\int_0^1 e^{\tilde{a}t} \cos(\tilde{b}t)\,dt = \frac{e^{\tilde{a}}(\tilde{a}\cos\tilde{b} + \tilde{b}\sin\tilde{b}) - \tilde{a}}{\tilde{a}^2 + \tilde{b}^2},$$

we calculate the integrals appearing in the rationality condition, and we call them $r_{11}$, $r_{12}$, $r_{21}$, and $r_{22}$:

$$\int_0^1 \xi_1(t)\omega_1(t)\,dt = \int_0^1 e^{(a+1)t}\sin(bt)\,dt = \begin{cases} r_{11}, & \text{if } (a+1)^2 + b^2 \neq 0, \\ 0, & \text{if } (a+1)^2 + b^2 = 0, \end{cases}$$

$$\int_0^1 \xi_1(t)\omega_2(t)\,dt = \begin{cases} r_{12}, & \text{if } (a+2)^2 + b^2 \neq 0, \\ 0, & \text{if } (a+2)^2 + b^2 = 0, \end{cases}$$

$$\int_0^1 \xi_2(t)\omega_1(t)\,dt = \begin{cases} r_{21}, & \text{if } (a+1)^2 + b^2 \neq 0, \\ 1, & \text{if } (a+1)^2 + b^2 = 0, \end{cases}$$

$$\int_0^1 \xi_2(t)\omega_2(t)\,dt = \begin{cases} r_{22}, & \text{if } (a+2)^2 + b^2 \neq 0, \\ 1, & \text{if } (a+2)^2 + b^2 = 0. \end{cases}$$

The computation of

$$\operatorname{sign} y(1) = \operatorname{sign}\left(\sum_{l=1}^{2} \alpha_l \sum_{j=1}^{2} g_j \int_0^1 \xi_l(t)\omega_j(t)\,dt\right)$$

is divided into three cases:

- Case $(a+1)^2 + b^2 \neq 0$, $(a+2)^2 + b^2 \neq 0$:

$$\operatorname{sign} y(1) = \operatorname{sign}(\alpha_1 g_1 r_{11} + \alpha_1 g_2 r_{12} + \alpha_2 g_1 r_{21} + \alpha_2 g_2 r_{22}) = \operatorname{sign}\left(\frac{P_1}{Q_1}\right).$$

- Case $(a+1)^2 + b^2 = 0$, $(a+2)^2 + b^2 \neq 0$:
$$\operatorname{sign} y(1) = \operatorname{sign}(\alpha_1 g_2 r_{12} + \alpha_2 g_1 + \alpha_2 g_2 r_{22}) = \operatorname{sign}\left(\frac{P_2}{Q_2}\right).$$

- Case $(a+1)^2 + b^2 \neq 0$, $(a+2)^2 + b^2 = 0$:
$$\operatorname{sign} y(1) = \operatorname{sign}(\alpha_1 g_1 r_{11} + \alpha_2 g_2 r_{21} + \alpha_2 g_2) = \operatorname{sign}\left(\frac{P_3}{Q_3}\right).$$

Thus we have three different expressions of the form $\frac{\mathbf{P}}{\mathbf{Q}}$.

Next we form the Boolean formula, $F = \operatorname{sign} y(1)$, evaluating polynomials $f_1 = (a+1)^2 + b^2 = 0$, $f_2 = (a+2)^2 + b^2 = 0$, $P_i > 0$, $Q_i > 0$, for $i \in \{1, 2, 3\}$. In the following table 1 means true and 0 means false for the above polynomial evaluation ($** = 1$ or $0$, i.e., extend the table).

| $f_1 = 0$ | $f_2 = 0$ | $P_1 > 0$ | $Q_1 > 0$ | $P_2 > 0$ | $Q_2 > 0$ | $P_3 > 0$ | $Q_3 > 0$ | $F$ |
|---|---|---|---|---|---|---|---|---|
| 0 | 0 | 1 | 1 | ** | ** | ** | ** | 1 |
| 0 | 0 | 1 | 0 | ** | ** | ** | ** | 0 |
| ⋮ | ⋮ | ⋮ | ⋮ | ⋮ | ⋮ | ⋮ | ⋮ | ⋮ |
| 1 | 0 | ** | ** | 1 | 1 | ** | ** | 1 |
| 1 | 0 | ** | ** | 1 | 0 | ** | ** | 0 |
| ⋮ | ⋮ | ⋮ | ⋮ | ⋮ | ⋮ | ⋮ | ⋮ | ⋮ |
| 0 | 1 | ** | ** | ** | ** | 0 | 0 | 1 |

In this case (see the statement of Goldberg-Jerrum bounds) $\Lambda = \{\alpha_1, \alpha_2, a, b, e^a \cos b, e^a \sin b\}$, $X = \{g_1, g_2\}$, $s = 8$, $d = 12$, and $l = 6$.

# B Functions Needed to Express an Exponential of a Real Matrix

Taking $2n^2$ functions $\xi_{r\ell}(t)$, $r = 1, \ldots, n$, $\ell = 1, \ldots, 2n$ to express $e^{At}$ in Proposition 5.1 is clear overcounting. In this appendix we estimate how many functions are needed to express $e^{At}$ for any $n \times n$ real matrix $A$.

Let us list the complex eigenvalues of $A$ as $a_1 \pm ib_1, \ldots, a_\kappa \pm ib_\kappa$, where $\kappa \leq \lfloor n/2 \rfloor$, in the decreasing order by their geometric multiplicities. We introduce functions

$$a_1 \pm ib_1: \quad \widetilde{M}_{a_1 \pm ib_1} := \{e^{a_1 t}\cos(b_1 t), te^{a_1 t}\cos(b_1 t), \ldots, t^{\lfloor \frac{n}{2}\rfloor - 1}e^{a_1 t}\cos(b_1 t),$$
$$e^{a_1 t}\sin(b_1 t), te^{a_1 t}\sin(b_1 t), \ldots, t^{\lfloor \frac{n}{2}\rfloor - 1}e^{a_1 t}\sin(b_1 t)\},$$

$$a_2 \pm ib_2: \quad \widetilde{M}_{a_2 \pm ib_2} := \{e^{a_2 t}\cos(b_2 t), \ldots, t^{\lfloor \frac{n}{4}\rfloor - 1}e^{a_2 t}\cos(b_2 t),$$
$$e^{a_2 t}\sin(b_2 t), \ldots, t^{\lfloor \frac{n}{4}\rfloor - 1}e^{a_2 t}\sin(b_2 t)\},$$

$$\vdots$$

$$a_\kappa \pm ib_\kappa: \quad \widetilde{M}_{a_\kappa \pm ib_\kappa} := \{e^{a_\kappa t}\cos(b_\kappa t),$$
$$e^{a_\kappa t}\sin(b_\kappa t)\}.$$

Linear combinations of the above functions can express all possible matrices $e^{At}$ when $A$ has purely complex eigenvalues.

As real eigenvalues can have higher geometric multiplicities, we add the following functions

$$M_{a_1 \pm ib_1} := \widetilde{M}_{a_1 \pm ib_1} \cup \{t^{\lfloor \frac{n}{2} \rfloor} e^{a_1 t}, \ldots, t^{n-1} e^{a_1 t}\},$$
$$M_{a_2 \pm ib_2} := \widetilde{M}_{a_2 \pm ib_2} \cup \{t^{\lfloor \frac{n}{4} \rfloor} e^{a_2 t}, \ldots, t^{\lfloor \frac{n}{2} \rfloor - 1} e^{a_2 t}\},$$
$$\vdots$$
$$M_{a_\kappa \pm ib_\kappa} := \widetilde{M}_{a_\kappa \pm ib_\kappa} \cup \{t e^{a_\kappa t}\}.$$

To express $e^{At}$ when $A$ has real eigenvalues $a_\ell$, $\ell = \kappa + 1, \ldots, n$ with geometric multiplicity one, we add at most $\lfloor n/2 \rfloor + 1$ functions $M_{a_\ell} = \{e^{a_\ell t}\}$ for $\ell = \kappa + 1, \ldots, n$.

The total number of functions in

$$M_{a_1 \pm ib_1} \cup \cdots \cup M_{a_\kappa \pm ib_\kappa} \cup M_{a_{\kappa+1}} \cup \cdots \cup M_{a_n}$$

is

$$2 \sum_{i=1}^{\lfloor n/2 \rfloor} \lfloor \frac{n}{2i} \rfloor + \sum_{i=1}^{\lfloor n/2 \rfloor} \left( \lfloor \frac{n}{i} \rfloor - \lfloor \frac{n}{2i} \rfloor \right) + \lfloor \frac{n}{2} \rfloor + 1$$
$$= \sum_{i=1}^{\lfloor n/2 \rfloor} \lfloor \frac{n}{2i} \rfloor + \sum_{i=1}^{\lfloor n/2 \rfloor} \lfloor \frac{n}{i} \rfloor + \lfloor \frac{n}{2} \rfloor + 1$$
$$\leq \frac{3n}{2} \sum_{i=1}^{\lfloor n/2 \rfloor} \frac{1}{i} + \lfloor \frac{n}{2} \rfloor + 1$$
$$\leq \frac{3n}{2} \left( \ln \lfloor \frac{n}{2} \rfloor + 1 \right) + \lfloor \frac{n}{2} \rfloor + 1 = O(n \ln n),$$

using the fact that $\ln(k+1) \leq \sum_{i=1}^{l} \frac{1}{i} \leq \ln(k) + 1$.

### Acknowledgments

This research has been financed in part by grants from the following organizations: the Academy of Finland, Finnish Academy of Sciences and Letters, Finnish Cultural Foundation and US Air Force (grant F47620-97-1-0159).